\documentclass{siamart251216}

\usepackage{amsfonts}
\usepackage{amssymb}
\usepackage{mathtools}
\usepackage{booktabs}
\usepackage{array}
\usepackage[most]{tcolorbox}
\usepackage{tikz}
\usepackage{enumitem}
\usepackage{overpic}

\newsiamremark{remark}{Remark}
\newsiamremark{workingnote}{Working note}
\newsiamthm{cert}{Certificate}
\newsiamthm{prop}{Proposition}

\newcommand{\R}{\mathbb{R}}
\newcommand{\C}{\mathbb{C}}

\definecolor{EmailFrame}{HTML}{6F675D}
\definecolor{EmailBack}{HTML}{FBFAF7}
\definecolor{EmailTitle}{HTML}{8A8176}
\definecolor{PromptFrame}{HTML}{2F5D85}
\definecolor{PromptBack}{HTML}{F4F8FC}
\definecolor{PromptTitle}{HTML}{386F9D}
\definecolor{ResponseFrame}{HTML}{3F735F}
\definecolor{ResponseBack}{HTML}{F4FAF7}
\definecolor{ResponseTitle}{HTML}{4B866F}

\newtcolorbox{emailbox}[1]{
  enhanced,
  breakable,
  colback=EmailBack,
  colframe=EmailFrame!65,
  coltitle=white,
  boxrule=0.4pt,
  borderline west={2pt}{0pt}{EmailTitle},
  arc=0.8mm,
  outer arc=0.8mm,
  left=8pt,
  right=8pt,
  top=10pt,
  bottom=7pt,
  before skip=7pt,
  after skip=8pt,
  fonttitle=\bfseries\itshape\small,
  attach boxed title to top left={xshift=7pt,yshift=-2mm},
  boxed title style={
    colback=EmailTitle,
    colframe=EmailTitle,
    boxrule=0pt,
    arc=0.8mm,
    outer arc=0.8mm,
    left=5pt,
    right=5pt,
    top=2pt,
    bottom=2pt},
  title={#1}}

\newtcolorbox{promptbox}[1]{
  enhanced,
  breakable,
  colback=PromptBack,
  colframe=PromptFrame!75,
  coltitle=white,
  boxrule=0.5pt,
  borderline east={2pt}{0pt}{PromptTitle},
  arc=1mm,
  outer arc=1mm,
  width=0.9\linewidth,
  enlarge left by=0.1\linewidth,
  left=7pt,
  right=8pt,
  top=11pt,
  bottom=7pt,
  before skip=8pt,
  after skip=6pt,
  fonttitle=\bfseries\small,
  attach boxed title to top right={xshift=-7pt,yshift=-2mm},
  boxed title style={
    colback=PromptTitle,
    colframe=PromptTitle,
    boxrule=0pt,
    arc=0.8mm,
    outer arc=0.8mm,
    left=5pt,
    right=5pt,
    top=2pt,
    bottom=2pt},
  title={#1}}

\newtcolorbox{responsebox}[1]{
  enhanced,
  breakable,
  colback=ResponseBack,
  colframe=ResponseFrame!75,
  coltitle=white,
  boxrule=0.5pt,
  borderline west={2pt}{0pt}{ResponseTitle},
  arc=1mm,
  outer arc=1mm,
  width=0.9\linewidth,
  enlarge right by=0.1\linewidth,
  left=8pt,
  right=7pt,
  top=11pt,
  bottom=7pt,
  before skip=3pt,
  after skip=8pt,
  fonttitle=\bfseries\small,
  attach boxed title to top left={xshift=7pt,yshift=-2mm},
  boxed title style={
    colback=ResponseTitle,
    colframe=ResponseTitle,
    boxrule=0pt,
    arc=0.8mm,
    outer arc=0.8mm,
    left=5pt,
    right=5pt,
    top=2pt,
    bottom=2pt},
  title={#1}}

\newtcolorbox{algorithmbox}[1]{
  enhanced,
  breakable,
  colback=black!2,
  colframe=black!45,
  coltitle=black,
  boxrule=0.45pt,
  borderline west={2pt}{0pt}{black!55},
  arc=0.8mm,
  outer arc=0.8mm,
  left=8pt,
  right=8pt,
  top=8pt,
  bottom=7pt,
  before skip=8pt,
  after skip=8pt,
  fonttitle=\bfseries\small,
  attach boxed title to top left={xshift=7pt,yshift=-2mm},
  boxed title style={
    colback=black!8,
    colframe=black!35,
    boxrule=0.35pt,
    arc=0.8mm,
    outer arc=0.8mm,
    left=5pt,
    right=5pt,
    top=2pt,
    bottom=2pt},
  title={#1}}

\makeatletter
\newcommand{\nopunctonce}{\gdef\@mypunct{\gdef\@mypunct{.}}}
\makeatother

\headers{Ten Digits on a Train}{M. J. Colbrook}
\title{Ten Digits on a Train:\\
AI-Assisted Verification of Two Eigenvalue Problems}
\author{Matthew J. Colbrook\thanks{Department of Applied Mathematics and Theoretical
Physics, University of Cambridge, Cambridge, CB3 0WA. (\email{mjc249@cam.ac.uk})}}

\begin{document}

\maketitle
\begin{center}
\small June 22, 2026\\[2pt]
\small\textit{Dedicated to Marco Marletta on his 60th birthday.}
\end{center}

\begin{abstract}
Accurate numerical eigenvalues are often difficult to certify, especially in singular or non-normal settings. This article reports a human--AI collaboration on two such computations. For a singular self-adjoint Schr\"odinger operator, a verified zero count and Dirichlet--Neumann bracketing certify the complete negative spectrum to ten decimal places. For a delicate non-normal atom--molecule benchmark, a previously unresolved resonance pair is separated, with each member enclosed to ten digits. The second result is achieved not by increasing the precision of one-way shooting, but by reformulating the problem as a global matching system for projective solution lines. The infinite tail is encoded as uncertainty in the terminal projective data, and a componentwise, tail-robust Krawczyk--Brouwer inclusion supplies the certificate. This gives a reusable architecture for analytic boundary-value systems with ill-conditioned propagation and uncertain asymptotic data. The collaboration also exposes the strengths and limits of AI assistance. AI rapidly produced accurate candidates and plausible proof strategies, but several failed, including one apparently complete tail argument that omitted the componentwise check required by a nonuniform polydisc. Validated computation is a stringent test of AI-assisted mathematics: the output is not merely a number, but a number with a proof. These examples show why the proof object matters, and why human mathematical judgment remained decisive. More broadly, as AI makes code, exposition, and plausible numerical claims inexpensive, standards for verification, attribution, peer review, and training must adapt. The implications are unsettling; the opportunity is extraordinary.
\end{abstract}



\section{Ten digits on a train}
\label{sec:intro}

Over the past few months, something has changed. Almost everywhere I visit, people are talking about the mathematical ability of AI systems. Sometimes the discussion is explicit; more often it sits just offstage. Which problems have become easier? Which parts of our work are being automated? How much of a serious computational project can now be initiated by a few prompts? This has become an elephant in the room for computational mathematics, and for mathematics more broadly. We talk about it in corridors, in offices, and over coffee, but not yet, I think, with the directness that the subject deserves.

The experiment reported here grew out of that atmosphere, but it began with a birthday-conference question. I asked Marco Marletta for a favorite eigenvalue, and then for a non-normal analog. The exchange was fun and fixed the mathematical starting point: a self-adjoint singular Coulomb problem and a resonance pair that he remembered as exceptionally difficult. Here are the four messages that set the scene.

\begin{emailbox}{Matthew to Marco}
\small\vspace{-1mm}
Hi Marco,
\par\smallskip
I hope that all is well! Quick question - what is your favourite eigenvalue?
\par\smallskip
Best, Matt\vspace{-2mm}
\end{emailbox}

\begin{emailbox}{Marco to Matthew}
\small\vspace{-1mm}
Hi Matt,
\par\smallskip
It is probably the ground state of hydrogen, so the bottom of the spectrum of the 1D Schroedinger operator with potential V(r) = 2/r\textasciicircum{}2 - 1/r. I believe its value is -0.0625, or -1/16, and I seem to remember my first computed approximation was -0.0621, before we had implemented the algorithm for doing better boundary conditions and meshing.
\par\smallskip
It was the first eigenvalue I ever computed numerically for a singular problem.
\par\smallskip
Many years later I was in Oxford visiting Nick Trefethen around the time he had just published Chebfun. He wanted to show me how good it was and asked me for a problem to test it. I gave him this problem, Chebfun promptly fell over, and Nick complained that it was not fair! Goodness, I had a lot to learn about social skills.
\par\smallskip
All the best, Marco\vspace{-2mm}
\end{emailbox}\vspace{-4mm}

\begin{emailbox}{Matthew to Marco}
\small\vspace{-1mm}
That is a wonderful story! Do you also have a favourite non-normal version of this?\vspace{-2mm}
\end{emailbox}\vspace{-4mm}

\begin{emailbox}{Marco to Matthew}
\small\vspace{-1mm}
I would choose one particular resonance pair from my 2014 many-author paper
\par\smallskip
https://orca.cardiff.ac.uk/id/eprint/64130/8/OA-20142015-45.pdf
\par\smallskip
It's the one at the bottom of page 15. It was an utter beast to compute. Never mind a 100 digit challenge, just getting 4 digits would have been nice!\vspace{-2mm}
\end{emailbox}

The first message leads to the $\ell=1$ sector of the radial Coulomb operator, with potential
$$\setlength\abovedisplayskip{6pt}\setlength\belowdisplayskip{6pt}
V_C(r)={2}/{r^2}-1/r,\qquad r>0,
$$
whose ground state is $-1/16=-0.0625$. Marco's recollection of an early numerical approximation near $-0.0621$, and of a later Chebfun \cite{DriscollHaleTrefethen2014} anecdote, made it a good control problem. The second message points to the resonance pair at the bottom of page 15 of B\"ogli, Brown, Marletta, Tretter, and Wagenhofer \cite{BoegliBrownMarlettaTretterWagenhoffer}, originally associated with Abramov, Aslanyan, and Davies \cite{AbramovAslanyanDavies}. Resonances encode oscillation and decay rates for open systems in quantum, wave, and scattering problems \cite{Zworski2017ScatteringResonances,dyatlov2019mathematical}.

I was on the train from Cambridge to Leeds when I decided to turn the above exchange into a test. \textit{Could an AI assistant be driven, by repeated mathematical questioning, to rigorous proof objects certifying ten digits of these eigenvalues?} By the end of the journey, I had a provisional answer. By the end of the day, I had what appeared to be a complete certificate. The AI system used in the computational exchanges was ChatGPT 5.5 (with Codex), with the intelligence setting ``Extra High'', in May 2026. I allowed it to search, write code, propose methods, make mistakes, and try to diagnose them. I did not allow a final claim to rest on the chat. Validated computation is a useful stress test for this reason. It asks the assistant not merely for a value, a program, or an explanation, but for a mathematically checkable object.

This paper reports that case study: not as a benchmark of an AI system in isolation, but as an account of a human--AI mathematical exchange in which numerical suggestions were repeatedly pushed toward independently checkable certificates. The examples should be read in that spirit. I analyze only the methods used by the large language model and the mathematical work needed to interpret, repair, or replace them. I do not claim that these are the best-conditioned, most economical, or most elegant validation strategies one could design from a blank page.

The paper sits at the intersection of two literatures. The first is the established tradition of validated computation: interval analysis, a posteriori estimates, Taylor models, and fixed-point tests have long turned numerical work into mathematical proof \cite{MooreKearfottCloud2009,Rump2010,Tucker,JohanssonArb}. The second literature is new and moving quickly: AI-assisted mathematics, including program-search discovery, human--AI proof search, and proof-assistant workflows whose statements and outputs must themselves be audited \cite{romera2024mathematical,collins2024evaluating,hubert2025olympiad}. The aim here is to examine their intersection in concrete validated spectral computations, with broader consequences for numerical mathematics.

The sharpest mathematical contribution is the resonance certificate. The calculation moves from the earlier two-decimal enclosure, which did not separate the two parity resonances, to two 17-digit certified enclosures. The key step is not extra precision. It is a change of proof object. Rather than validate a long non-normal shoot or a scalar determinant obtained from one-way propagation, the method projectivizes the solution line at each mesh point, treats all node coordinates as simultaneous unknowns, and verifies the resulting sparse global matching system by a Krawczyk inclusion.\footnote{The Krawczyk method \cite{krawczyk1969newton,Neumaier1990} combines a Newton-like correction with Jacobian interval bounds.} The tail estimate is then inserted as a set-valued final-row perturbation, and a robust Krawczyk--Brouwer inclusion gives the corresponding half-line resonance.\footnote{This should be distinguished from the validation in \cite{BoegliBrownMarlettaTretterWagenhoffer}, where interval arithmetic and the argument principle are applied to a finite-interval determinant.} This gives a reusable validation framework for analytic boundary-value systems with ill-conditioned propagation and uncertain terminal data.

\begin{figure}[t]
\centering
\begin{tikzpicture}[
  x=1cm,y=1cm,
  >=stealth,
  font=\small,
  every node/.style={align=center}
]
\tikzset{
  basebox/.style={
    rounded corners=1mm,
    minimum width=2.15cm,
    minimum height=0.68cm,
    inner xsep=4pt,
    inner ysep=2.5pt,
    line width=0.45pt
  },
  explorebox/.style={basebox, draw=PromptFrame!78, fill=PromptBack},
  assessbox/.style={basebox, draw=ResponseFrame!78, fill=ResponseBack},
  certbox/.style={basebox, draw=black!62, fill=black!3},
  lane/.style={
    rounded corners=1.4mm,
    line width=0.45pt
  },
  flowarrow/.style={
    ->,
    line width=1.05pt,
    draw=#1,
    line cap=round,
    shorten >=2pt,
    shorten <=2pt
  },
  arrowlabel/.style={
    font=\footnotesize,
    fill=white,
    inner xsep=2pt,
    inner ysep=1pt
  }
}

\draw[lane, draw=PromptFrame!45, fill=PromptBack!55]
  (-0.10,0.60) rectangle (10.65,2.05);
\node[font=\footnotesize\bfseries,text=PromptFrame!85!black,anchor=west]
  at (-0.10,2.30) {AI-assisted exploration};

\node[explorebox] (prompt) at (1.10,1.36) {prompt};
\node[explorebox] (reply) at (3.70,1.36) {response\\code/method};
\node[explorebox] (cand) at (6.30,1.36) {candidate\\or failure};
\node[assessbox] (assess) at (9.45,1.36) {author\\assessment};

\draw[flowarrow=PromptFrame!88!black] (prompt) -- (reply);
\draw[flowarrow=PromptFrame!88!black] (reply) -- (cand);
\draw[flowarrow=PromptFrame!88!black] (cand) -- (assess);
\draw[flowarrow=PromptFrame!88!black]
  (assess.south) .. controls (8.55,0.72) and (2.05,0.72) ..
  node[arrowlabel,text=PromptFrame!88!black,pos=0.55,below=1pt] {new prompt}
  (prompt.south);

\draw[lane, draw=black!36, fill=black!2]
  (-0.10,-1.82) rectangle (10.65,-0.42);
\node[font=\footnotesize\bfseries,text=black!72,anchor=east]
  at (10.50,-0.18) {mathematical certification};

\node[certbox] (cert) at (2.38,-1.14) {finite\\certificate};
\node[certbox] (test) at (5.28,-1.14) {theorem or\\validation test};
\node[certbox] (claim) at (8.18,-1.14) {mathematical\\claim};

\draw[flowarrow=black!68] (cert) -- (test);
\draw[flowarrow=black!68] (test) -- (claim);

\coordinate (certjoin) at (6.30,0.24);

\draw[flowarrow=black!68]
  (cand.south) -- (certjoin)
  -- node[arrowlabel,text=black!65,midway,below=1pt] {if checkable}
  (2.38,0.24)
  -- (cert.north);

\draw[flowarrow=black!68]
  (assess.south) -- (9.45,0.24) -- (certjoin);
	
\end{tikzpicture}\vspace{-2mm}
\caption{The workflow used in the two examples. The assistant operates in an exploration loop, where prompts lead to proposed methods, code, numerical candidates, and useful failures. A computed spectral claim enters the paper only after passing through a separate certification layer.\vspace{-3mm}}
\label{fig:workflow}
\end{figure}
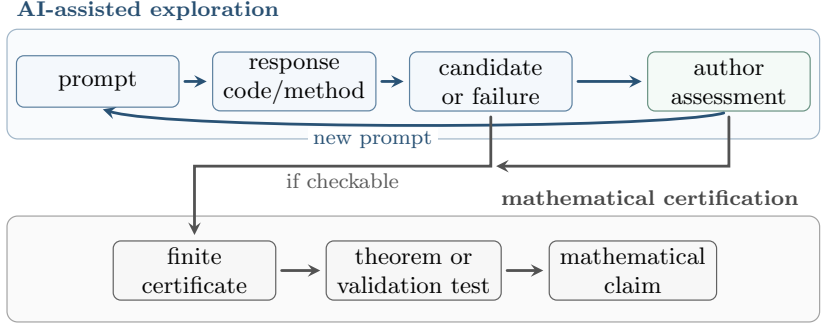

Throughout, I show a working style for interacting with an AI assistant summarized in \Cref{fig:workflow}. The assistant operates in the exploration loop, where prompts produce methods, code, candidates, and failures. The mathematical claims pass through a separate certification layer. In the resonance problem, this separation was essential: failed approaches led to the global projective formulation, and a draft certificate that looked finished failed an audit because it controlled global maxima where the proof required nodewise inclusion in a polydisc. AI can and does make mistakes. The lesson is not to dismiss AI-assisted work, but to make every claim checkable.

The prompt and response\footnote{I also record the response time to roughly indicate the relative difficulty of each stage. These timings are machine-, server load-, and internet-speed-dependent (UK train Wi-Fi is notoriously bad) and are longer than the run times of the final interval-arithmetic certificate codes.} boxes below record this exchange. The explanations in the paper reconstruct and correct the intervening work from the assistant's files---scripts, notes, output tables, and certificates---and assess which mathematical claims they support. I have checked all of the code for the certificates line by line.

\begin{table}[t]
\caption{Numerical certificates verified to ten digits. The displayed values are rounded to ten decimal places. The resonance validation below has an error radius less than $7\times 10^{-19}$. The self-adjoint operator considered is a modification of the exactly solvable Coulomb operator.\vspace{-1mm}}
\label{tab:headline}
\centering
\begingroup
\footnotesize
\setlength{\tabcolsep}{5pt}
\newcommand{\headlineDash}{
  \multicolumn{3}{@{}c@{}}{
    \tikz[baseline=-0.6ex]
      \draw[black!36,dash pattern=on 3.5pt off 2.4pt,line width=0.35pt]
        (0,0) -- (0.985\linewidth,0);}\\[-2pt]
}
\begin{tabular*}{\linewidth}{@{\extracolsep{\fill}}lcl@{}}
\toprule
certified object & displayed value & verification method \\
\midrule
modified self-adjoint ground state
&$-0.0531967615$
&
Temple certificate \\
\headlineDash
modified self-adjoint negative spectrum
&
11 intervals in \Cref{tab:selfadjoint-spectrum}
&
Sturm verifier \\
\headlineDash
odd/Dirichlet resonance
&
$0.9180018174-17.0006882497\,\mathrm{i}$
&
projective Krawczyk \\
\headlineDash
even/Neumann resonance
&
$0.9178768454-17.0006448370\,\mathrm{i}$
&
projective Krawczyk \\
\bottomrule
\end{tabular*}
\endgroup\vspace{-3mm}
\end{table}

\Cref{tab:headline} summarizes the numerical outcomes. The article is also intended as a guided exposition of validated numerics: the certificates are introduced through the questions they answer, rather than as a catalog of interval algorithms. The three certificates become progressively more complicated: a variational Temple certificate for the ground state, a global Sturm and bracketing certificate for the full negative spectrum, and a non-variational projective Krawczyk certificate for the resonance pair. I learned some interesting mathematics during the exchange, which I hope the reader will enjoy as well.

\textit{Where does this leave us?} Wider implications are discussed in \Cref{sec:comparison_final,sec:conclusion}, but the numerical analysis community should not treat this as a curiosity at the edge of the field. The point is not to argue against AI: the opportunities are extraordinary, and the exchange below let me explore methods, failures, and reformulations at a speed that would have been hard to imagine only recently. Precisely because the opportunity is real, the standards matter. AI changes the cost of trying ideas, the apparent difficulty of examples, and the ease with which plausible numerical claims can be produced. If our standards for validation, credit, and research training do not adapt, polished but weakly checked claims may travel faster than understanding.

The structure of the paper is as follows. \Cref{sec:coulomb,sec:selfadjoint} treat the self-adjoint problems: first, a control example, and then a harder certificate problem. \Cref{sec:resonance} gives the non-normal resonance problem, the failed validation routes, and the global projective Krawczyk certificate. \Cref{sec:comparison_final} discusses lessons from the certificates, including the social consequences of AI-assisted computational mathematics. We conclude in \Cref{sec:conclusion}.

\section{A control problem: Marletta's Coulomb eigenvalue}
\label{sec:coulomb}

The first problem Marletta gave was the reduced radial problem on $(0,\infty)$:
\begin{equation}\setlength\abovedisplayskip{6pt}\setlength\belowdisplayskip{6pt}
\label{easy_problem}
    -u''(r)+\left(\frac{2}{r^2}-\frac1r\right)u(r)=\lambda u(r),
    \qquad u(r)=O( r^2)\quad(r\downarrow 0).
\end{equation}
This is the $\ell=1$ Coulomb problem. The Coulomb spectrum of hydrogen was one of the proving grounds of quantum mechanics. Balmer's empirical formula and Rydberg's generalization organized the observed lines into integer-labeled series \cite{Balmer1885,Rydberg1890}; Bohr's model explained the same pattern by quantized orbits \cite{Bohr1913}; and Schr\"odinger recast hydrogen as an eigenvalue problem and recovered the spectrum from a wave equation \cite{Schrodinger1926}. After separation of variables, the angular momentum $\ell=1$ case has centrifugal term $\ell(\ell+1)/r^2=2/r^2$, giving exactly \eqref{easy_problem} in the present units.

The essential spectrum of this operator is $[0,\infty)$ \cite[Thm.~10.9]{Teschl}. To find the negative eigenvalues, let $\lambda=-\kappa^2$, with $\kappa>0$, and set $x=2\kappa r$. If
$
    u(r)=r^2 e^{-\kappa r} y(x),
$
then the differential equation becomes
$$\setlength\abovedisplayskip{6pt}\setlength\belowdisplayskip{6pt}
    x y''+(4-x)y'
      +\left(\frac{1}{2\kappa}-2\right)y=0.
$$
This is Kummer's equation \cite[Sec.~13.2(i)]{OlverLozierBoisvertClark2010}. The solution regular at $r=0$ is square integrable at infinity only when the hypergeometric series terminates, that is when
$$\setlength\abovedisplayskip{6pt}\setlength\belowdisplayskip{6pt}
    \frac{1}{2\kappa}-2=m,\qquad m=0,1,2,\ldots .
$$
Then $y$ is a generalized Laguerre polynomial $L_m^{(3)}(x)$. Setting $n=m+1\in\mathbb{N}$ gives infinitely many negative eigenvalues accumulating at the bottom of the essential spectrum:
$$\setlength\abovedisplayskip{6pt}\setlength\belowdisplayskip{6pt}
    \lambda_n=-\kappa^2=-\frac{1}{4(n+1)^2},\qquad n=1,2,\ldots .
$$
Since $x=2\kappa r=r/(n+1)$, the corresponding eigenfunctions are, up to normalization,
\begin{equation}\setlength\abovedisplayskip{3pt}\setlength\belowdisplayskip{6pt}
\label{hydrogen_efuns}
    u_n(r)=r^2 e^{-r/(2(n+1))}L_{n-1}^{(3)}(r/(n+1)),
    \qquad n=1,2,\ldots .
\end{equation}
Thus the smallest eigenvalue is $-1/16$, with eigenfunction proportional to $r^2e^{-r/4}$.

A short, earlier AI exchange, not reproduced here, recovered this answer quickly, but not by merely printing a remembered formula. I asked for a certified numerical computation. For a positive chosen trial function $v$, the assistant found that the pointwise quotient
$
    \smash{\frac{-v''+(2/r^2-1/r)v}{v}}
$
was identically $-1/16$. When I next asked for enclosures of the ten lowest eigenvalues, the attempted numerical validation exposed the terminating Laguerre structure. This relegated the Coulomb problem to a useful control case. Once the right ansatz is found, the proof reduces to polynomial identities. The more revealing test was to modify the problem in a new conversation that did not contain this derivation.

\section{A modified self-adjoint problem}
\label{sec:selfadjoint}

\subsection{The lowest eigenvalue}

To remove the exact Laguerre structure while keeping the same singular half-line setting, I sent the following prompt.

\begin{promptbox}{Prompt 1: modifying Marletta's eigenvalue}
\small
\setlength{\parindent}{1.5em}

Below, I will paste an email exchange with Marco Marletta. He describes his favourite eigenvalue. Let's change the potential to $\frac{2}{r^2}-\frac{1}{r(1+r/200)^2}$. I want you to investigate whether this has a negative eigenvalue. If so, I want you to come up with a way to compute its smallest eigenvalue. You will need to rigorously prove error bounds on the eigenvalue and fully describe your computational procedure. You will find Python installed in the current directory. You can probably build a rational computation verification in Python, or you may need to look at interval arithmetic implementations.
\end{promptbox}

For self-adjoint operators, and here in particular, a good approximate eigenvector can often be turned into a rigorous eigenvalue bound by combining residual estimates, comparison theorems, and variational principles (see \Cref{cert:selfadjoint} below). 

The modified operator in this prompt is
$$\setlength\abovedisplayskip{6pt}\setlength\belowdisplayskip{6pt}
    H=-\frac{\mathrm{d}^2}{\mathrm{d}r^2}+V(r),\qquad V(r)=\frac{2}{r^2}
      -\frac{1}{r(1+r/200)^2}
$$
on $L^2(0,\infty)$. The modification preserves the singular radial character of Marletta's problem, but it removes the exact Laguerre solution. Near $0$, the dominant equation remains $-u''+\frac{2}{r^2}u=0$, whose power-law solutions are $r^2$ and $r^{-1}$. The admissible branch therefore satisfies $u(r)=O(r^2)$, while the $r^{-1}$ branch is excluded. As for the Coulomb operator, the essential spectrum of $H$ is $[0,\infty)$ \cite[Thm.~10.9]{Teschl}. If a trial function has Rayleigh quotient below $0$, then the min--max principle implies the existence of a discrete negative eigenvalue.

\begin{responsebox}{AI response, \textit{worked for 10m 54s}}
\small\setlength{\parindent}{1.5em}
Yes. The modified potential still has a negative ground-state eigenvalue.
I produced a rigorous certificate in \path{EIGENVALUE.md} and the verifier script
\path{verify_temple.py}. Certified enclosure:
$$\setlength\abovedisplayskip{6pt}\setlength\belowdisplayskip{6pt}
 -0.05319676154302054294
 \leq \lambda_1 \leq
 -0.05319676154233948813.
$$
The proof uses a five-term hydrogenic Ritz vector, computes its Rayleigh quotient and residual with rational interval arithmetic, then applies Temple's inequality. The spectral gap needed for Temple comes from comparing the new potential with Marco's Coulomb potential, whose second $l=1$ level is $-1/36$. I also ran the verifier successfully with the local Python.
\end{responsebox}

At this stage these were claims made by the assistant, not yet certified results. I have audited the proposed argument and verifier, and this is the same pattern used throughout. Here, the assistant first used a finite-difference discretization and then formed a five-dimensional Ritz approximation that found a candidate near $-0.0532$. It then proposed a certificate: round the Ritz vector to rational data, compute a Rayleigh quotient and residual bound, use the Coulomb operator as a comparison, and apply Temple's inequality with an explicit gap.

The certificate is small enough to describe explicitly. It uses the trial vector
\begin{equation}\setlength\abovedisplayskip{6pt}\setlength\belowdisplayskip{6pt}
\label{ground_state_trial}
    u(r)=r^2 e^{-r/4}+\sum_{j=1}^4 c_j r^2 e^{-r/4}(r/10)^j.
\end{equation}
The transcript does not establish why this form was selected. Mathematically, it is a natural polynomial--exponential enrichment of the first Coulomb eigenfunction $u_1$ in \cref{hydrogen_efuns}. The rational coefficients were obtained by rounding an exploratory Ritz vector:
$$\setlength\abovedisplayskip{6pt}\setlength\belowdisplayskip{6pt}
\begin{array}{rclcrcl}
c_1&=&\displaystyle
 \frac{911357530019777039}{50000000000000000000000},&&
c_2&=&\displaystyle
 \frac{3464720637087477951}{500000000000000000000},\\[2mm]
c_3&=&\displaystyle
 -\frac{653458902923720157}{2500000000000000000000},&&
c_4&=&\displaystyle
 \frac{1159101098623161979}{50000000000000000000000}.
\end{array}
$$
The displayed rational coefficients define the trial function exactly; they need not enclose the floating-point Ritz coefficients. The exploratory calculation is used only to select a good trial function. The floating-point Ritz calculation has no further role.

For the normalized vector $v=u/\|u\|$, the verifier computes enclosures for the Rayleigh quotient $\rho=\langle Hv,v\rangle$ and squared residual $\delta^2=\|(H-\rho)v\|^2$. The functions in the expansion \cref{ground_state_trial} are products of powers of $r$ and $e^{-r/4}$.
After integration by parts and the substitution $t=r+200$, the only non-elementary quantity needed is the single scalar
$$\setlength\abovedisplayskip{6pt}\setlength\belowdisplayskip{6pt}
    C=e^{100}E_1(100)=\int_0^\infty \frac{e^{-s}}{100+s}\,\mathrm{d}s.
$$
For example, the inner products involve
$$\setlength\abovedisplayskip{6pt}\setlength\belowdisplayskip{6pt}
    K=\int_0^\infty \frac{e^{-r/2}r^3}{(r+200)^2}\,\mathrm{d}r =e^{100}\int_{200}^{\infty}e^{-t/2}
       \left(t-600+\frac{120000}{t}-\frac{8000000}{t^2}\right)\,\mathrm{d}t .
$$
The nonnegative powers of $t$ are evaluated by repeated integration by parts, while
$$\setlength\abovedisplayskip{6pt}\setlength\belowdisplayskip{6pt}
    e^{100}\int_{200}^{\infty}\frac{e^{-t/2}}{t}\,\mathrm{d}t=C,
    \qquad
    e^{100}\int_{200}^{\infty}\frac{e^{-t/2}}{t^2}\,\mathrm{d}t
      =\frac1{200}-\frac{C}{2}.
$$
Thus, $K=-40796+4120000\,C.$ Similarly, the verifier represents each required integral as $a+bC$, with $a,b\in\mathbb Q$. The scalar $C$ is enclosed by the finite identity
$$\setlength\abovedisplayskip{6pt}\setlength\belowdisplayskip{6pt}
\frac{1}{100+s}
=
 \sum_{k=0}^{N-1}\frac{(-1)^k s^k}{100^{k+1}}
 +\frac{(-1)^N s^N}{100^N(100+s)} .
$$
After multiplying by $e^{-s}$ and integrating, this gives
$$\setlength\abovedisplayskip{6pt}\setlength\belowdisplayskip{6pt}
    C=\sum_{k=0}^{N-1}\frac{(-1)^k k!}{100^{k+1}}+R_N,
    \qquad
    0\leq (-1)^N R_N\leq \frac{N!}{100^{N+1}} .
$$
This is a rigorously bounded asymptotic expansion: the terms initially decrease, but factorial growth eventually dominates as $N\to\infty$. Nevertheless, with $N=60$, the verifier obtains an interval of width less than
$8.33\times10^{-41}$, which proves that
$$
    C=0.009901942286733018406405931819804200422\ldots .
$$
For example, $K=-40796+4120000\,C=0.0022213\ldots$ with large cancellation, which illustrates why exact rational endpoints are preferable here.

The reduction to one constant is a useful example of certificate design. It is more transparent, smaller, and easier to check than asking an interval quadrature routine to integrate each required integral independently. The variational structure then provides a global enclosure.

\begin{cert}[Ground-state enclosure]
\label{cert:selfadjoint}
For the modified Coulomb operator $H$, the lowest eigenvalue $\lambda_1$ satisfies
$$\setlength\abovedisplayskip{6pt}\setlength\belowdisplayskip{6pt}
-0.05319676154302054294
\leq \lambda_1 \leq
-0.05319676154233948813.
$$
\end{cert}

\begin{proof}
The verifier rigorously gives $\delta^2\leq 1.7311721\times10^{-14}$ and
$$\setlength\abovedisplayskip{6pt}\setlength\belowdisplayskip{6pt}
\rho\in[-0.0531967615423394881331,\,
        -0.0531967615423394881330].
$$
Given a self-adjoint operator $\mathcal{L}$ that is bounded below, let $\lambda_1(\mathcal{L})\leq\lambda_2(\mathcal{L})\leq\cdots$ denote its discrete eigenvalues below the essential spectrum. If there are a finite number $N$ of these, set $\lambda_{N+j}(\mathcal{L})$ to be the bottom of the essential spectrum for $j\in\mathbb{N}$. These are the min--max levels of $\mathcal{L}$. Since $-{1}/(r(1+r/200)^2)\geq-1/r$, comparison with the $\ell=1$ Coulomb operator $H_C$ in \Cref{sec:coulomb} implies that $\lambda_2(H)\geq\lambda_2(H_C)=-1/36.$ Rayleigh--Ritz gives $\lambda_1(H)\leq \rho<-1/36$. Therefore $H$ has exactly one eigenvalue below $-1/36$ and Temple's inequality \cite[Thm.~1]{temple1928theory} gives
$$\setlength\abovedisplayskip{6pt}\setlength\belowdisplayskip{6pt}
    \rho-\frac{\delta^2}{-1/36-\rho}\leq \lambda_1(H)\leq \rho.
$$
Substituting the outward-rounded intervals gives the displayed enclosure.
\end{proof}

In the above proof, the condition $\rho<\lambda_2(H_C)\leq \lambda_2(H)$ lets a single trial vector certify the bottom eigenvalue, rather than merely locating a nearby spectral point.

\nopunctonce
\subsection{Finite or infinite negative spectrum?}
\label{sec:finite_or_infinite_evals}

At this stage, I was very excited that this could all be achieved at the push of a button. (It took me significantly longer to write the above explanation of the certificate and to check it.) In the usual mathematical way, I immediately wanted more.

\begin{promptbox}{Prompt 2: finite or infinite negative spectrum}
\small\setlength{\parindent}{1.5em}
Wonderful! Are there infinitely many negative eigenvalues or finitely many? If there are finitely many, say $N$, compute the first $\min(10,N)$ negative eigenvalues to 10 digits. If there are infinitely many, compute the smallest 10 to 10 digits. Again, everything has to be rigorously justified.
\end{promptbox}

\begin{responsebox}{AI response, \textit{worked for 21m 56s}}
\small\setlength{\parindent}{1.5em}
It is finite: there are exactly $11$ negative eigenvalues. I wrote the justification in \path{SPECTRUM.md}, with the shooting generator \path{spectrum_shooting.py}. The first ten, ordered from most negative upward, are:
$-0.0531967615$,
$-0.0193674097$,
$-0.0082420874$,
$-0.0036685443$,
$-0.0016120468$,
$-0.0006697369$,
$-0.0002515677$,
$-0.0000804935$,
$-0.0000198995$, and
$-0.0000031084$.
The proof structure is: the potential is nonnegative for $r\geq9800+4000\sqrt6$, so there can only be finitely many negative eigenvalues; the zero-energy regular solution has exactly $11$ zeros, giving $N=11$ by Sturm oscillation; Dirichlet--Neumann bracketing then encloses the first ten half-line eigenvalues in intervals of radius $5\times10^{-12}$.
\end{responsebox}

The assistant's first observation was that $V(r)\geq0$ for $r\geq r_+=9800+4000\sqrt6$. In fact, since
$\lim_{r\to\infty}r^2\,V(r)= 2>-1/4,
$
there is a Kneser-type half-line criterion \cite[Cor.~3.6(ii)]{BehrndtSchmitzTeschlTrunk2023} that rules out an infinite sequence of negative eigenvalues accumulating at $0$. Hence, there are only finitely many negative eigenvalues. The second observation involves $u_0$, the regular solution of the zero-energy equation
\begin{equation}\setlength\abovedisplayskip{6pt}\setlength\belowdisplayskip{6pt}
\label{zero_energy_equation}
    -u_0''(r)+\left(\frac{2}{r^2}
      -\frac{1}{r(1+r/200)^2}\right)u_0(r)=0,
    \qquad u_0(r)\sim r^2 \quad(r\downarrow 0).
\end{equation}
Sturm oscillation identifies the number of negative eigenvalues of $H$ with the number of zeros of $u_0$, provided that zero is not itself an eigenvalue \cite[Thm.~2.1, Cor.~2.4, and the discussion following Thm.~7.4]{gesztesy1996zeros}.

However, I should highlight a \textbf{warning}. The initial response overclaimed since the produced code \path{spectrum_shooting.py} used floating-point arithmetic. In preparing this paper, I separated the exploratory computation from the proof step. The certified claims below use the later Arb \cite{JohanssonArb} verifier \path{selfadjoint_interval_verify.py}, which checks the zero crossings, the tail signs, and the Dirichlet--Neumann endpoint signs by interval arithmetic. It proves one zero in each of the following intervals:
$$\setlength\abovedisplayskip{6pt}\setlength\belowdisplayskip{6pt}
\begin{array}{rclcrcl}
1&:&[10.681579,10.683579],&
2&:&[26.241254,26.243254],\\
3&:&[49.959273,49.961273],&
4&:&[85.068026,85.070026],\\
5&:&[137.168771,137.170771],&
6&:&[216.494771,216.496771],\\
7&:&[343.195032,343.197032],&
8&:&[561.492538,561.494538],\\
9&:&[985.487724,985.489724],&
10&:&[1999.515565,1999.517565],\\
11&:&[5906.061779,5906.063779].
\end{array}
$$
This is a validated shooting calculation. The verifier starts at $\varepsilon=10^{-2}$ from a regular Frobenius expansion with an explicit bound on the omitted tail, and a coefficient bound shows that there is no zero on $(0,\varepsilon)$. The verifier then propagates the interval enclosure for $(u_0,u_0')$ by Taylor steps for
$
    u_0''=Vu_0,
$
bounding each remainder by a Cauchy--Picard majorant. A detailed description is delayed until \Cref{sec:eval_enclosures}.

The verifier first proves a sign change across each of the 11 displayed intervals. Independently, it covers $[\varepsilon,T]$ by certified Taylor steps. On each step, either the range enclosure of $u_0$ excludes zero, or the endpoint signs are opposite and the range enclosure of $u_0'$ excludes zero. The scan therefore finds exactly 11 zeros and no others. Since the displayed boxes are disjoint and each contains a zero, each contains exactly one zero. It is verified that $u_0(T)<0$ and $u_0'(T)<0$. Since $V(r)\geq0$ for $r\geq T$ and $u_0''=Vu_0\leq0$ for $u_0<0$, this sign information rules out further zeros on $[T,\infty)$. \Cref{fig:u0-plots} (right) illustrates this behavior. This also proves that zero is not an eigenvalue: $\lim_{r\to\infty}u_0(r)=-\infty$ so $u_0$ cannot be square integrable.

\subsection{Analytic form and zero counting missed by the AI}

The zero-energy solution $u_0$ has a closed form and analytic argument for its 11 zeros, which the assistant did not notice. This does not invalidate the interval certificate above; rather, it is a useful reminder that, just like a human, an AI-assisted route can find a checkable proof without finding the simplest analytic description that often helps a proof.

Set
$
    r=200\tan^2\!\left({\theta}/{2}\right)$ and
$u_0(r)=\sqrt{r}\,y(\theta)$ for $0<\theta<\pi$.
A direct calculation transforms the zero-energy equation into
\begin{equation}\setlength\abovedisplayskip{6pt}\setlength\belowdisplayskip{6pt}
    -\frac{1}{\sin\theta}
      \frac{\mathrm d}{\mathrm d\theta}
      \left(\sin\theta\,y'(\theta)\right)
    +\frac{9}{\sin^2\theta}y(\theta)
    =200\,y(\theta).
\label{eq:associated-legendre}
\end{equation}
This is the associated Legendre equation of order $3$ and degree
$$\setlength\abovedisplayskip{6pt}\setlength\belowdisplayskip{6pt}
    \nu={(\sqrt{801}-1)}/{2},
    \qquad
    \nu(\nu+1)=200
$$
\cite[Secs.~14.2--14.3]{OlverLozierBoisvertClark2010}.
The solution with the required behavior $u_0(r)\sim r^2$ as
$r\downarrow0$ is
$$\setlength\abovedisplayskip{6pt}\setlength\belowdisplayskip{6pt}
    u_0(r)
    =
    6\times 200^{3/2}\sqrt r\,
    \mathsf P_{\nu}^{-3}
    \left(\frac{200-r}{200+r}\right),
$$
where $\mathsf P_{\nu}^{-3}$ denotes the Ferrers function on
$(-1,1)$.

This transformation also gives an analytic zero count. The Sturm--Liouville operator on the left-hand side of
\eqref{eq:associated-legendre} has eigenvalues $n(n+1)$, $n=3,4,5,\ldots$.
Since
$
    13\times14=182<200<210=14\times15,
$
Sturm oscillation shows that the solution regular at $\theta=0$ has exactly $11$ zeros in $(0,\pi)$. Because $r=200\tan^2(\theta/2)$ maps $(0,\pi)$ bijectively onto $(0,\infty)$, the zero-energy solution $u_0$ has exactly $11$ positive zeros.

\begin{figure}[t]
\centering
\includegraphics[width=0.5\linewidth]{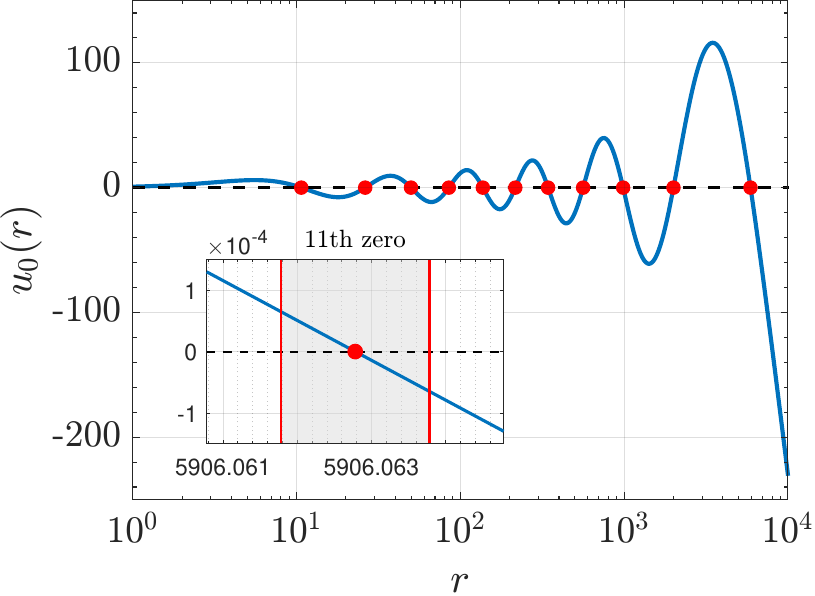}\hfill
\includegraphics[width=0.49\linewidth]{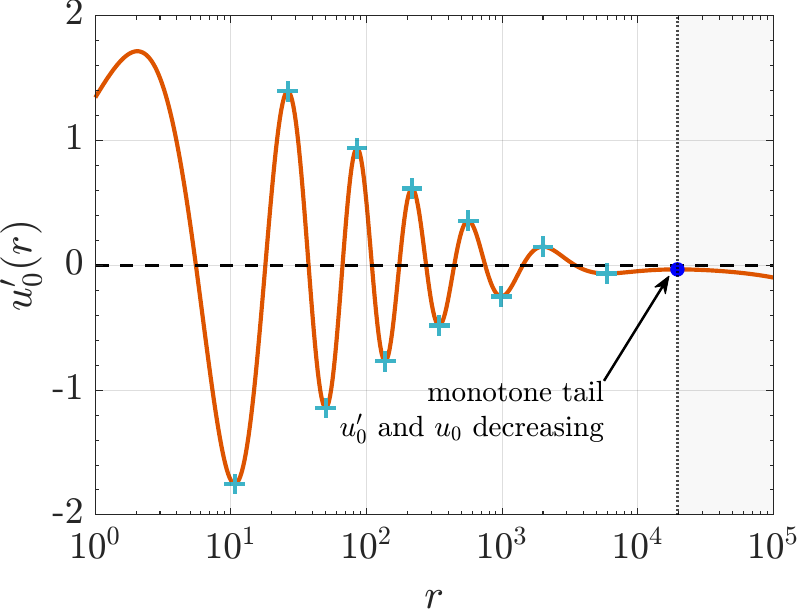}\vspace{-2mm}
\caption{The regular solution of \eqref{zero_energy_equation}. Left: $u_0$, with the verified interval for the eleventh zero shaded in the inset. Right: $u_0'$, with the beginning of the monotone tail indicated. The zeros of $u_0$ are shown as red dots in the left panel and crosses in the right panel.\vspace{-4mm}}
\label{fig:u0-plots}
\end{figure}

\Cref{fig:u0-plots} shows this solution and its derivative. The left panel plots $u_0$, with a magnified panel near the eleventh zero. The right panel plots $u_0'$ and marks the onset of the nonnegative tail. For $r>0$, since $u_0''=V u_0$, the zeros of $u_0$, together with the zeros $r_\pm=9800\pm4000\sqrt6$ of the potential, are the critical points of $u_0'$. At the certified point $T>r_+$, the verified signs of $u_0$ and $u_0'$ give the non-oscillatory tail used above.

\subsection{Enclosures for all negative eigenvalues}
\label{sec:eval_enclosures}

Before giving the final argument for the eigenvalue enclosures, I will record my final prompt for this problem---there was only one sensible thing left to do.\footnote{The reader may find the tone of some prompts amusing. I often find (and am often told) that AI assistants work better when encouraged.}\vspace{-1mm}

\begin{promptbox}{Prompt 3: the eleventh eigenvalue}
\small\setlength{\parindent}{1.5em}
Wow! Great job. Since there are 11, we may as well compute the 11th as well. Please can you do this. Rigorously argue why your computations are correct, explicitly describe the method and any computations you have
done.
\end{promptbox}

\begin{responsebox}{AI response, \textit{worked for 5m 29s}}
\small\setlength{\parindent}{1.5em}
Done. The eleventh eigenvalue is certified in the interval
$$\setlength\abovedisplayskip{6pt}\setlength\belowdisplayskip{6pt}
[-1.6240198497433\times10^{-7},\,
  -1.6239198497433\times10^{-7}].
$$
I updated the spectrum note with all eleven eigenvalues in \path{SPECTRUM.md}. This includes the finiteness argument, Sturm oscillation count, Dirichlet/Neumann bracketing, Frobenius initialization, and interval-shooting certification.
\end{responsebox}

As in the preceding response, this was an overclaim at the time: \path{spectrum_shooting.py} supplied a plausible candidate interval, but the certification was established only by the later Arb verifier \path{selfadjoint_interval_verify.py}.

The task is to enclose eigenvalues of $H$ between the finite-interval Neumann and Dirichlet eigenvalues \cite[Sec.~2.3 and Sec.~4.4]{Teschl}. Let $H_L^{\rm D}$ and $H_L^{\rm N}$ be the finite-interval realizations on $(0,L)$, with the same $u(r)=O(r^2)$ condition at $0$ and with Dirichlet and Neumann conditions respectively at $L$.\footnote{The floating-point locator \path{spectrum_shooting.py} used the larger truncation point $25000$ to find candidate intervals. The rigorous interval-arithmetic script \path{selfadjoint_interval_verify.py} then verifies the required endpoint signs at
$
    L=T=19597.9589711328>r_+=9800+4000\sqrt6,
$
which already lies in the nonnegative tail of the potential.} Let $H_{L,\infty}^{\rm D}$ and $H_{L,\infty}^{\rm N}$ denote the tail realizations on $(L,\infty)$. The corresponding quadratic-form domains satisfy
$
\mathcal Q(H_L^{\rm D}\oplus H_{L,\infty}^{\rm D})
\subset
\mathcal Q(H)
\subset
\mathcal Q(H_L^{\rm N}\oplus H_{L,\infty}^{\rm N}),
$
and the forms agree on their common domains. Consequently,
$
H_L^{\rm N}\oplus H_{L,\infty}^{\rm N}
\leq H\leq
H_L^{\rm D}\oplus H_{L,\infty}^{\rm D}
$
in the form sense. Since $V(r)\geq0$ for $r\geq L$, both tail operators are nonnegative. The min--max principle therefore gives
\begin{equation}\setlength\abovedisplayskip{6pt}\setlength\belowdisplayskip{6pt}
  \lambda_j(H_L^{\rm N})
  \leq\lambda_j(H)
  \leq\lambda_j(H_L^{\rm D}),
  \qquad j=1,\ldots,11.
\label{bracket_L}
\end{equation}

The verifier uses a Frobenius expansion to construct rigorous initial balls for $u(\varepsilon,\lambda)$ and $u'(\varepsilon,\lambda)$, where $\varepsilon=10^{-2}$. Near $r=0$,
$$\setlength\abovedisplayskip{6pt}\setlength\belowdisplayskip{6pt}
    \frac{1}{r(1+r/200)^2}
    =
    \sum_{m=0}^{\infty}
    \frac{(-1)^m(m+1)}{200^m}r^{m-1}.
$$
Writing
$u(r,\lambda)=r^2\sum_{n=0}^{\infty}a_n(\lambda)r^n$ with $a_0=1$
and comparing the coefficient of $r^n$ in the eigenvalue equation gives
$$\setlength\abovedisplayskip{6pt}\setlength\belowdisplayskip{6pt}
    n(n+3)a_n+
    \sum_{m+j=n-1}\frac{(-1)^m(m+1)}{200^m}a_j
    +\lambda a_{n-2}=0,
$$
with the convention that $a_k=0$ for $k<0$. For every spectral endpoint used below, $|\lambda|<0.054$. The verifier initializes the regular solution at $\varepsilon=10^{-2}$ by retaining the coefficients $a_0(\lambda),\ldots,a_{80}(\lambda)$. The recurrence above implies, by induction, that
$$\setlength\abovedisplayskip{6pt}\setlength\belowdisplayskip{6pt}
    |a_n(\lambda)|
    \leq
		\frac{|\lambda|}{n(n+3)}+
    \frac{1}{n(n+3)}
        \sum_{m=0}^{n-1}\frac{m+1}{200^m} 
    \leq
    \frac{{40000}/{39601}+0.06}{4}
    <1,
    \qquad n\geq1.
$$
Consequently, the omitted tails beginning with $a_{81}$ satisfy
$$\setlength\abovedisplayskip{6pt}\setlength\belowdisplayskip{6pt}
\left|
u(\varepsilon,\lambda)
{-}\sum_{n=0}^{80}a_n(\lambda)\varepsilon^{n+2}
\right|
{\leq}
\frac{\varepsilon^{83}}{1-\varepsilon},\quad
\left|
u'(\varepsilon,\lambda)
{-}\sum_{n=0}^{80}(n+2)a_n(\lambda)\varepsilon^{n+1}
\right|
{\leq}\frac{83\varepsilon^{82}-82\varepsilon^{83}}{(1-\varepsilon)^2}.
$$
The verifier evaluates the two finite sums in Arb ball arithmetic, adds these symmetric tail bounds, and uses the resulting balls as the initial data at $r=\varepsilon$.

The verifier then propagates the solution by validated local Taylor expansions. At a point $x$, let
$$\setlength\abovedisplayskip{6pt}\setlength\belowdisplayskip{6pt}
    q(x+t,\lambda)=V(x+t)-\lambda
       =\sum_{k=0}^{\infty}q_k t^k,
    \qquad
    u(x+t,\lambda)=\sum_{n=0}^{\infty}b_n t^n.
		$$
The equation $u''=qu$ gives the recurrence
$$\setlength\abovedisplayskip{6pt}\setlength\belowdisplayskip{6pt}
    b_{n+2}
    =
    \frac{1}{(n+1)(n+2)}
    \sum_{k=0}^{n}q_k b_{n-k},
    \qquad n\geq0.
$$
To bound the omitted Taylor tail over the interval $0\leq t\leq h$, choose $\rho>h$ such that the disk $|z-x|\leq\rho$ avoids the singularity at $0$, and set
$$\setlength\abovedisplayskip{6pt}\setlength\belowdisplayskip{6pt}
    Q=
    \frac{2}{(x-\rho)^2}
    +\frac{40000}{(x-\rho)(x+200-\rho)^2}
    +|\lambda|.
$$
If $Q\rho^2/2<1$, the integral equation and a Picard majorant give
$$\setlength\abovedisplayskip{6pt}\setlength\belowdisplayskip{6pt}
    \sup_{|z-x|\leq\rho}|u(z)|
    \leq
    M:=
    \frac{U+\rho U_1}{1-Q\rho^2/2},
$$
where $U$ and $U_1$ are rigorous upper bounds for $|u(x)|$ and $|u'(x)|$. Cauchy's estimate then yields, with $\tau=h/\rho$,
$$\setlength\abovedisplayskip{6pt}\setlength\belowdisplayskip{6pt}
\left|\sum_{n=d+1}^{\infty}b_n h^n\right|
\leq
M\frac{\tau^{d+1}}{1-\tau},\qquad
\left|\sum_{n=d+1}^{\infty}n b_n h^{n-1}\right|
\leq
\frac{M}{\rho}
\frac{\tau^d\bigl(d+1-d\tau\bigr)}{(1-\tau)^2}.
$$
The verifier therefore obtains rigorous enclosures of the endpoint values and of the ranges of $u$ and $u'$ throughout every step.

\begin{figure}[t]
\centering
\includegraphics[width=0.49\linewidth]{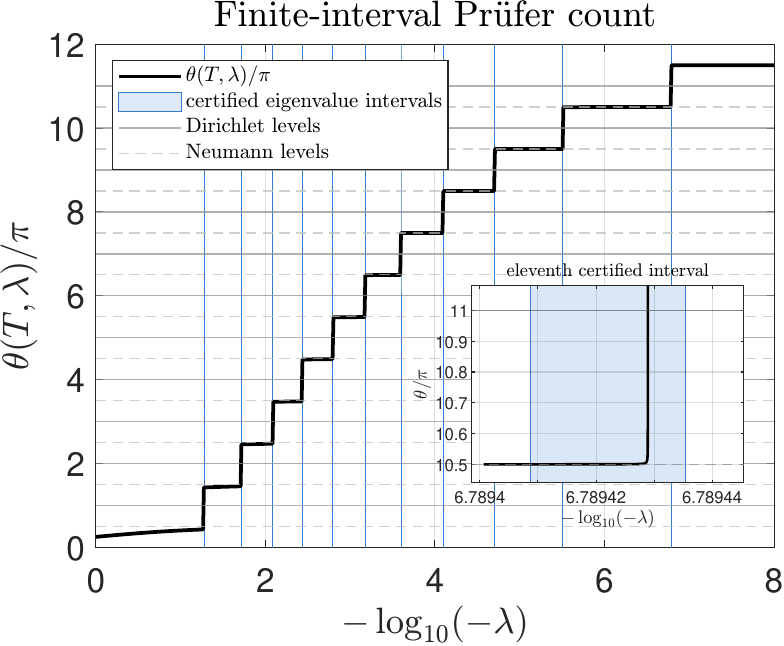}\hfill
\begin{overpic}[width=0.46\linewidth]{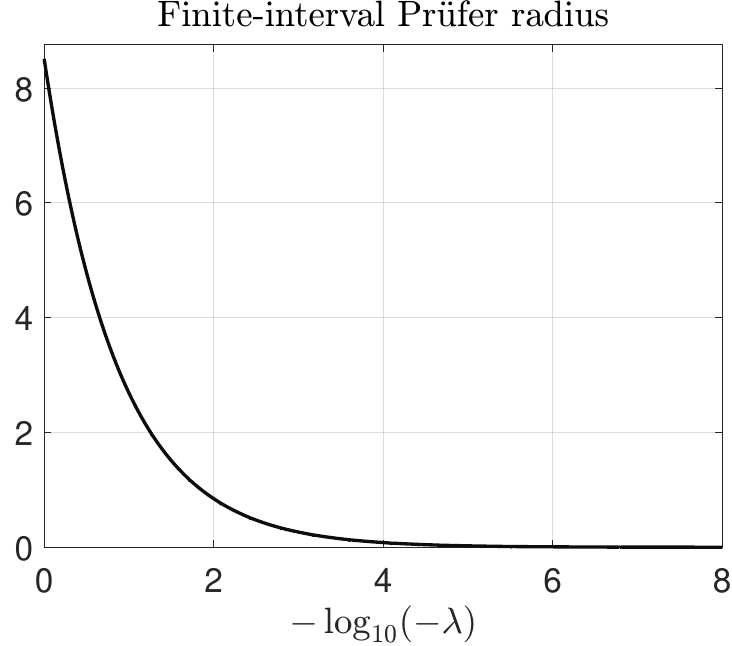}
  \put(-6,25){\rotatebox{90}{\small $10^{-3}\log_{10}\mathfrak R(T,\lambda)$}}
\end{overpic}\vspace{-2mm}
\caption{Floating-point illustration of the finite-interval Pr\"ufer
calculation at $T=19597.9589711328$. Left: $\theta(T,\lambda)/\pi$ as a
function of $-\log_{10}(-\lambda)$. Integer levels correspond to
Dirichlet conditions $u(T,\lambda)=0$, and half-integer levels correspond
to Neumann conditions $u'(T,\lambda)=0$. The blue bands are the certified
common bracketing intervals $m_j\pm5\times10^{-12}$; their endpoints are
not the crossing locations. The inset magnifies the eleventh interval and
marks the approximate Dirichlet and Neumann crossings obtained from the
floating-point Pr\"ufer data. Right: $10^{-3}\log_{10}\mathfrak R(T,\lambda)$,
where $\mathfrak R^2=u(T,\lambda)^2+u'(T,\lambda)^2$. Its broad growth is
expected: for $\lambda=-\kappa^2<0$, the tail equation is approximately
$u''=\kappa^2u$, so a generic terminal solution grows like $e^{\kappa T}$.
The figure is illustrative only; the certificate
uses interval endpoint signs and Sturm oscillation.}
\label{fig:prufer-count}
\end{figure}

The finite-interval Dirichlet and Neumann eigenvalues are the zeros of $u(L,\lambda)$ and $u'(L,\lambda)$, respectively. To count them, define the Pr\"ufer radius and angle by
$$\setlength\abovedisplayskip{6pt}\setlength\belowdisplayskip{6pt}
u(r,\lambda)=\mathfrak{R}(r,\lambda)\sin\theta(r,\lambda),
\qquad
u'(r,\lambda)=\mathfrak{R}(r,\lambda)\cos\theta(r,\lambda),
$$
with $\theta(r,\lambda)\to0$ as $r\downarrow0$. A Wronskian calculation gives
$$\setlength\abovedisplayskip{6pt}\setlength\belowdisplayskip{6pt}
\frac{\partial\theta}{\partial\lambda}(L,\lambda)=
\frac{\int_0^L u(r,\lambda)^2\,\mathrm{d}r}
{u(L,\lambda)^2+u'(L,\lambda)^2}>
0.
$$
Moreover, $\theta(L,\lambda)\to0$ as $\lambda\to-\infty$. Since $u_0(T)<0$ and $u_0'(T)<0$ for $L=T=19597.9589711328$, and $u_0$ has 11 positive zeros,
$$\setlength\abovedisplayskip{6pt}\setlength\belowdisplayskip{6pt}
11\pi<\theta(T,0)<{23\pi}/{2}.
$$
Hence, the conditions $\theta(T,\lambda)\in\pi\mathbb Z$ give exactly 11 negative Dirichlet eigenvalues, and the conditions $\theta(T,\lambda)\in\pi/2+\pi\mathbb Z$ give exactly 11 negative Neumann eigenvalues. \Cref{fig:prufer-count} illustrates the Sturm count: $\theta(r,\lambda)$ records the endpoint boundary conditions, while $\mathfrak{R}(r,\lambda)$ records the exponential growth of a generic terminal solution.

The code \path{spectrum_shooting.py} supplied approximate roots from which the candidate intervals were selected. The Arb verifier \path{selfadjoint_interval_verify.py} repeats the shooting calculation with interval Taylor arithmetic. For each interval $m_j\pm5\times10^{-12}$, it propagates the regular solution at both $\lambda$-endpoints to evaluate $u(T,\lambda)$ and $u'(T,\lambda)$. Opposite signs of $u(T,\lambda)$ at the two endpoints imply a Dirichlet finite-interval eigenvalue in the interval, and opposite signs of $u'(T,\lambda)$ imply a Neumann one. Since the verifier finds an eigenvalue in each of 11 pairwise disjoint negative intervals, every interval contains exactly one Dirichlet and one Neumann eigenvalue. The inequality \cref{bracket_L} then encloses $\lambda_j(H)$. The end result is the following certificate.

\begin{table}[t]
\caption{All negative eigenvalues of $H$. The Arb verifier certifies the Sturm count and the Dirichlet--Neumann endpoint signs for the intervals $m_j\pm 5\times10^{-12}$, with row $j$ corresponding to $\lambda_{j}$. The first eigenvalue also has the independent Temple certificate in \Cref{cert:selfadjoint}.}
\label{tab:selfadjoint-spectrum}
\centering\footnotesize
\begin{tabular}{@{}rlll@{}}
\toprule
$j$ & midpoint of certified interval & radius & value to 10 decimal places \\
\midrule
1  & $-5.3196761542354\times10^{-2}$ & $5\times10^{-12}$ & $-0.0531967615$ \\
2  & $-1.9367409728912\times10^{-2}$ & $5\times10^{-12}$ & $-0.0193674097$ \\
3  & $-8.2420874282006\times10^{-3}$ & $5\times10^{-12}$ & $-0.0082420874$ \\
4  & $-3.6685443298679\times10^{-3}$ & $5\times10^{-12}$ & $-0.0036685443$ \\
5  & $-1.6120468156717\times10^{-3}$ & $5\times10^{-12}$ & $-0.0016120468$ \\
6  & $-6.6973685120550\times10^{-4}$ & $5\times10^{-12}$ & $-0.0006697369$ \\
7  & $-2.5156770074481\times10^{-4}$ & $5\times10^{-12}$ & $-0.0002515677$ \\
8  & $-8.0493482670423\times10^{-5}$ & $5\times10^{-12}$ & $-0.0000804935$ \\
9  & $-1.9899482902437\times10^{-5}$ & $5\times10^{-12}$ & $-0.0000198995$ \\
10 & $-3.1084439935000\times10^{-6}$ & $5\times10^{-12}$ & $-0.0000031084$ \\
11 & $-1.6239698497433\times10^{-7}$ & $5\times10^{-12}$ & $-0.0000001624$ \\
\bottomrule
\end{tabular}
\end{table}

\begin{cert}[Negative spectrum]
\label{prop:selfadjoint-eleven}
The modified half-line operator has exactly eleven negative eigenvalues contained in the intervals in \Cref{tab:selfadjoint-spectrum}.
\end{cert}

\section{The non-normal resonance pair}
\label{sec:resonance}

The second problem was deliberately less forgiving. The target is the resonance pair $\hat\lambda_2^\pm$ from \cite[Table~6]{BoegliBrownMarlettaTretterWagenhoffer}. That work already proved the existence of one resonance in each parity sector, but computed both to two digits and did not separate them. Our aim is to replace that common enclosure by two disjoint enclosures determining both components to ten decimal places.

Here, the feedback loop in \Cref{fig:workflow} became even more important: high-precision candidates appeared quickly, but the failed certificates forced changes in the proof strategy. The result is a general and reusable method: a global projective multiple-shooting certificate that turns a non-normal complex-scaled resonance calculation from an unstable shooting problem into a verified finite-dimensional existence theorem.

\subsection{Background}

The underlying operator is the Schr\"odinger operator
$$\setlength\abovedisplayskip{6pt}\setlength\belowdisplayskip{6pt}
   -\frac{\mathrm{d}^2}{\mathrm{d}x^2}+Q(x),\qquad
    Q(x)=(x^2-1.6)\exp(-x^2/10),\qquad x\in\R.
$$
However, the numbers in \cite[Table~6]{BoegliBrownMarlettaTretterWagenhoffer} are not $L^2$ eigenvalues, but resonances.

Because $Q$ is an even function, the full-line problem decomposes into two half-line problems: odd states satisfy the Dirichlet condition $u(0)=0$, while even states satisfy the Neumann condition $u'(0)=0$. With $\sqrt{-\lambda}$ chosen so that $\operatorname{Re}\sqrt{-\lambda}>0$, the condition at infinity is exponential decay, $u(x)\sim \exp(-\sqrt{-\lambda}x)$, which gives the usual $L^2$ eigenvalue problem. Resonances arise as poles of the meromorphic continuation across the continuous spectrum $\lambda\in[0,\infty)$, equivalently as discrete eigenvalues of a complex-scaled operator discussed below. After continuation, the same outgoing solution need not lie in $L^2$; in the examples below, it grows exponentially on the real axis. This is the complex-scaling interpretation used in \cite[Sec.~2, eqs.~(2.1)--(2.5)]{BoegliBrownMarlettaTretterWagenhoffer}; see also the classical references \cite{AguilarCombes1971,SimonComplexScaling1978} and the book \cite{dyatlov2019mathematical}. The complex-scaled boundary-value problems that we use below are highly non-normal, so even plausible numerical roots computed with high precision can be treacherous.

\begin{promptbox}{Prompt 4: from Marco's desired four digits to ten}
\small\setlength{\parindent}{1.5em}
Wonderful, now let's do a different problem. Marco sent me the message pasted below. I want us to push the boat out. Figure out a way to compute this resonance to 10 digits of accuracy (not just the 4 digits that Marco would like). It has to be fully rigorous and accurate. We may need to come up with new ideas to do this.
\end{promptbox}

\begin{responsebox}{AI response, \textit{worked for 1h 15m 12s}}
\small\setlength{\parindent}{1.5em}
I've identified the target and obtained high-precision candidate values. The pair is the Abramov pair $\hat\lambda_2^-,\hat\lambda_2^+$ in Table 6, at the bottom of page 15 of B\"ogli--Brown--Marletta--Tretter--Wagenhofer. Using the paper's stretched scaling $(R,\theta,X)=(10,0.76,10)$, I get:
\begin{gather}\setlength\abovedisplayskip{6pt}\setlength\belowdisplayskip{6pt}
\hat\lambda_2^- \approx
0.9180018173874698-17.00068824973865\,\mathrm{i},\\
\hat\lambda_2^+ \approx
0.9178768454495876-17.00064483695551\,\mathrm{i}.
\end{gather}
These are high-precision candidates, not yet a rigorous certificate. The route I think will get us the rigorous 10 digits is: work in the $z$-plane, shoot backwards from $X=10$ with the outgoing condition $y'(X)=-\sqrt{-z}$, form the Dirichlet and Neumann determinants, propagate the ODE and variational equation with interval Taylor arithmetic, apply an interval Newton/Krawczyk test on small boxes in $z$, use the paper's truncation bound $\alpha_{X,\theta,R}\leq1.77\times10^{-17}$, so that the tail error is far below the requested ten-digit scale. The next step is to turn the Taylor shooting into an interval certificate.
\end{responsebox}

Here, $\hat\lambda_2^-$ and $\hat\lambda_2^+$ correspond to the Dirichlet and Neumann half-line problems, respectively. The assistant identified both resonances and suggested the following possible validation route. Ordinary complex scaling sends $\lambda$ to $z=e^{2\mathrm{i}\theta}\lambda$. For this difficult pair, one instead chooses $R>0$ and follows the complex ray $\xi=R e^{\mathrm{i}\theta}x$ for $x\geq0$, with $y(x)=u(\xi)$. Multiplying the transformed equation by $R^2e^{2\mathrm{i}\theta}$ and truncating the domain gives
$$\setlength\abovedisplayskip{6pt}\setlength\belowdisplayskip{6pt}
    -y''(x)+q_{\theta,R}(x)y(x)=z\,y(x),\qquad 0\leq x\leq X,
$$
where
$$\setlength\abovedisplayskip{6pt}\setlength\belowdisplayskip{6pt}
 q_{\theta,R}(x)
 =R^2e^{2\mathrm{i}\theta}
   \bigl(R^2e^{2\mathrm{i}\theta}x^2-1.6\bigr)
   \exp\!\left(-\frac{R^2e^{2\mathrm{i}\theta}x^2}{10}\right)\quad\text{and}\quad
	z=R^{2}e^{2\mathrm{i}\theta}\lambda.
$$
We use the same parameters as \cite{BoegliBrownMarlettaTretterWagenhoffer}: $X=10$, $\theta=0.76\,\mathrm{rad}$, and $R=10$.

\begin{figure}[t]
\centering
\includegraphics[width=0.49\linewidth]{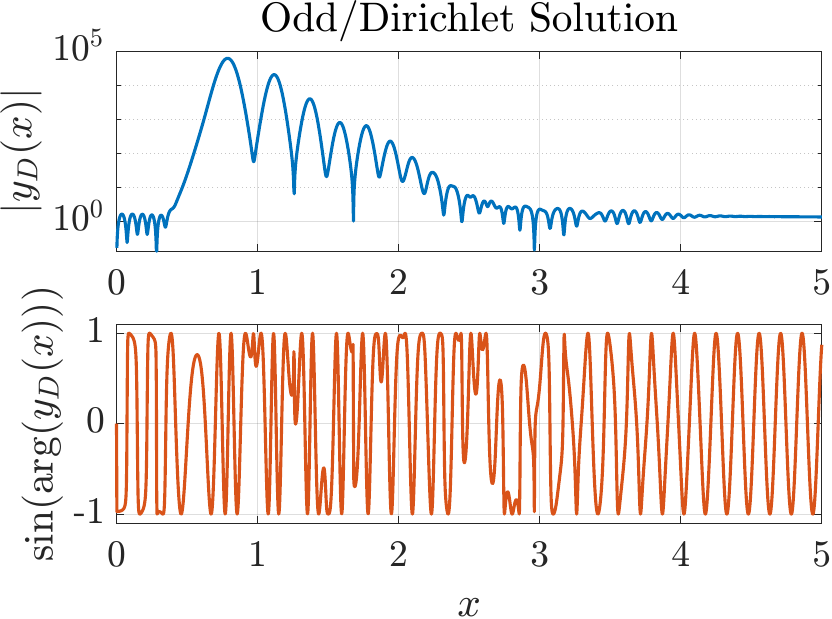}
\includegraphics[width=0.49\linewidth]{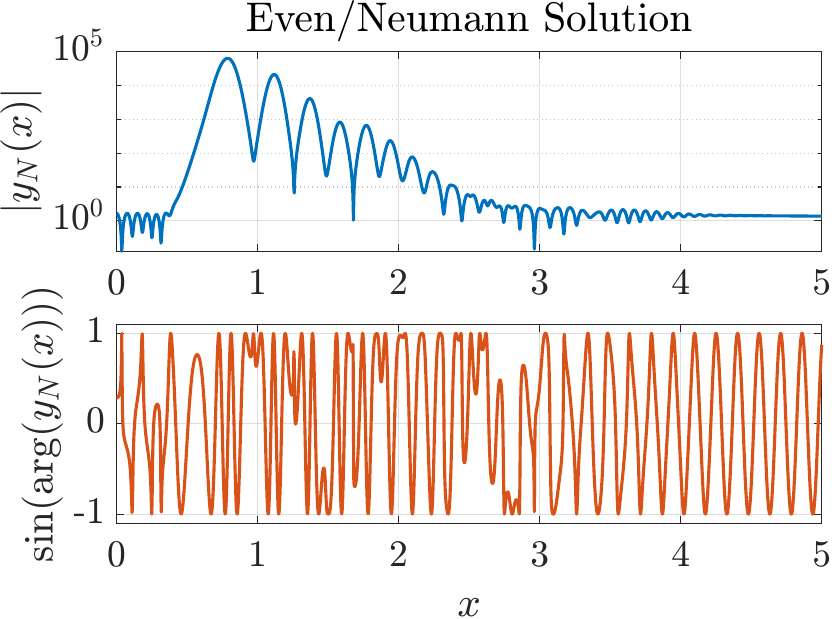}\vspace{-4mm}
\caption{Finite-interval resonance solutions for the stretched complex-scaled problem. The solutions are normalized by $y_D(10)=1$ and $y_N(10)=1$, and only the interval $0< x\leq5$ is shown to make the oscillation and growth visible. The upper panels plot $|y(x)|$ on a logarithmic scale; the lower panels plot $\sin(\arg y(x))$. This figure is illustrative only: the rigorous resonance enclosure is supplied by the global projective Krawczyk certificate.\vspace{-4mm}}
\label{fig:resonance_solutions}
\end{figure}

Let $k(z)=\sqrt{-z}$, with $\operatorname{Re} k(z)>0$. The outgoing free solution in the absence of a potential is the decaying exponential $\exp(-k(z)\,x)$. Thus, we impose the following boundary condition at the truncation point $X$:
\begin{equation}\setlength\abovedisplayskip{6pt}\setlength\belowdisplayskip{6pt}
\label{tau_BC}
y'(X)+k(z)y(X)=0.
\end{equation}
For each $z$, integrating the outgoing solution backwards from $X$ to $0$ yields $F_D(z)=y(0;z)$ and $F_N(z)=y'(0;z)$. Zeros of $F_D$ and $F_N$ give the odd and even \textit{finite-interval} resonances, respectively. \Cref{fig:resonance_solutions} shows the corresponding solutions at the numerical candidate resonances. A tail bound can then be used to account for the contribution of the tail of $q_{\theta,R}$ on $[X,\infty)$ to the true outgoing solution.

In exact arithmetic, this is a clean analytic formulation. For computation, this is an unstable one-way shooting: terminal and local interval errors are propagated through the non-normal dynamics to $x=0$.

\subsection{Useful failures}

The final certificate did not appear at once. Several failures are summarized in \Cref{tab:failures}. Each failed attempt found an obstruction and changed the next question (\Cref{fig:workflow}). They also had a clear logic. Amplitudes were too unstable. Changes of basis and midpoint matching improved the scaling, but still left a boundary determinant built from long propagation. Projective coordinates revealed the right variable (the solution line), but not the right proof. The final method in \Cref{resonance_success} kept the projective variable and replaced one-way shooting with a global multiple-shooting interval Newton system.

The AI assistant repeatedly predicted that the next modification would complete the proof. The subsequent failures show why such predictions were overly optimistic.

\begin{table}[t]
\caption{Diagnostic routes to the resonance certificate. Each failed route identified which object should not be certified or an obstacle to overcome, and so narrowed the form of the final proof.\vspace{-1mm}}
\label{tab:failures}
\centering
\begingroup
\footnotesize
\setlength{\tabcolsep}{2pt}
\newcommand{\headlineDash}{
  \multicolumn{3}{@{}c@{}}{\vspace{-4pt}
    \tikz[baseline=-1.8ex]
      \draw[black!36,dash pattern=on 3.5pt off 2.4pt,line width=0.35pt]
        (0,0) -- (0.985\linewidth,0);}\\[-3.5pt]
}
\begin{tabular}{@{}>{\raggedright\arraybackslash}p{0.27\linewidth}
                >{\raggedright\arraybackslash}p{0.36\linewidth}
                >{\raggedright\arraybackslash}p{0.34\linewidth}@{}}
\toprule
route & obstruction & lesson \\
\midrule
amplitude shooting &
enclosures grew rapidly in the unstable solution component &
raw amplitudes are poor variables for long validated propagation \\

\headlineDash

variation-of-constants and amplitude Lohner forms &
free growth and some wrapping were reduced, but one-way
amplification remained &
better local coordinates did not cure the long composition \\

\headlineDash

midpoint or Wronskian matching &
the enclosure of the matching determinant remained too wide &
this characteristic formulation was still too ill-conditioned \\

\headlineDash

Riccati/projective shooting &
amplitude growth disappeared, but a chart denominator could no
longer be certified to exclude zero &
solution lines are the appropriate variables, but the long one-way
composition remained unsuitable \\

\headlineDash
centered projective Lohner forms &
local projective remainders initially stayed small but later widened &
the node lines should be treated as simultaneous unknowns \\

\headlineDash

more precision &
roundoff shrank while the principal widths and failure point persisted &
the tested formulation was not roundoff limited \\

\midrule

global projective multiple shooting &
short validated maps coupled by projective matching equations &
\textbf{a componentwise Krawczyk inclusion succeeded}\\
\bottomrule
\end{tabular}
\endgroup\vspace{-4mm}
\end{table}

\begin{promptbox}{Prompt 5: asking for interval arithmetic}
\small\setlength{\parindent}{1.5em}
This is excellent work! Now let's push for a fully rigorous verification. Will you need to do interval arithmetic in Python?
\end{promptbox}

\begin{responsebox}{AI response, \textit{instant response}}
\small\setlength{\parindent}{1.5em}
Yes, I think we do need interval arithmetic in Python, or something equivalent. For this problem, ordinary high-precision floating point is only enough to locate the candidates. To make a 10-digit claim rigorous, we need certified enclosures for: the ODE propagation, the determinant values $F_\pm(z)$, the derivative with respect to $z$, the interval Newton or Krawczyk inclusion test, and the truncation/scaling error from replacing the half-line problem by the $X=10$ boundary problem.

The good news is that the problem is low-dimensional enough that we probably do not need a huge interval package. I'd build a small complex-interval Taylor solver in Python using rational/decimal endpoints where safe, with our own interval complex square-root branch for $\sqrt{-z}$. The key validation would be:
$
N(\mathcal B)=z_0 - F'(\mathcal B)^{-1}F(z_0)
$
or a Krawczyk variant, showing $N(\mathcal B)\subset \mathcal B$ for each small box $\mathcal B$ around the two candidate roots.
\end{responsebox}

On each short interval, the Taylor coefficients give a finite, checkable enclosure of the transfer map and its derivative. The first issue was whether these local enclosures could be composed across $[0,X]$ without losing the boundary condition.

\begin{promptbox}{Prompt 6: committing to the interval route}
\small\setlength{\parindent}{1.5em}
Wonderful, let's do it!
\end{promptbox}

\begin{responsebox}{AI response, \textit{worked for 1h 27m 14s}}
\small\setlength{\parindent}{1.5em}
Yes, and I've started the rigorous verification machinery. I installed python-flint and built several Arb-based verifier prototypes in the workspace.

Important: I do not yet have a completed rigorous 10-digit certificate. The naive interval methods fail: direct amplitude propagation creates huge wrapping/unstable-direction enclosures, even though the floating point determinant is tiny. The projective Riccati version is the promising route, but its current chart switching still needs a sharper Lohner/Taylor-model style enclosure; otherwise Arb refuses reciprocal chart switches when the interval still contains zero. That is a genuine rigor issue, not a coding nuisance.

So the next step is clear: replace the current absolute interval Riccati propagation by a centered projective Taylor model. In practice, propagate a high-precision midpoint line plus a small rigorous error ball in a locally chosen M\"obius chart, updating the chart center along the computed trajectory. That should remove the zero-containing chart-switch problem and give a usable Krawczyk proof.
\end{responsebox}

The first attempts used direct interval shooting, propagating the solution in complex balls using interval Taylor methods with rigorous remainder estimates at each step. The difficulty was that these local errors compounded: they were fed into transverse directions amplified by the non-normal shooting dynamics. The resulting enclosures, therefore, measured not only uncertainty in the resonance but also the instability of the proof architecture itself: small Taylor and interval errors carried through an unstable one-way propagation.

Variation of constants removed much of the free evolution, while interior matching shortened the shooting intervals. Both improved the scaling, but the argument still ended with a scalar matching determinant built from interval-dependent data, whose enclosure was too wide to certify the zero. 

Projective coordinates removed the irrelevant amplitude and followed the line spanned by $(y,y')$, the object on which the boundary conditions depend. However, in the one-way Riccati and projective attempts, the intervals became too wide to prove that the denominators in maps induced by the transfer matrices stay away from zero.

At this stage, it was natural to ask whether more precision could solve the issue.

\begin{promptbox}{Prompt 7: extended precision as a possible cure?}
\small\setlength{\parindent}{1.5em}
OK, would it help to use extended precision floats? Can we do it by just using way more precision? If not, try this other method you discussed.
\end{promptbox}

\begin{responsebox}{AI response, \textit{worked for 40m 21s}}
\small\setlength{\parindent}{1.5em}
Extended precision alone won't help. I tested exactly that on the centered-projective verifier: same failure point at 360 bits and 900 bits. The rounding part shrank dramatically, but the actual enclosing ball did not. At the failure point: (a) 360 bits; projective coordinate radius about 21.1635544598, denominator radius about 85.0013815165 (b) 900 bits; projective coordinate radius still about 21.1635544598, denominator radius still about 85.0013815165. So precision improved the tiny Arb roundoff decorations from around $10^{-105}$ to $10^{-267}$, but the useful interval width stayed unchanged. That means the obstruction is interval dependency/wrapping, not floating point precision.

I tried the other method and added a new verifier that: builds midpoint frames along the left and right solution lines; converts each rigorous Taylor step into a local projective M\"obius transfer map; propagates only the centered projective coordinate and its $z$-derivative; tries a Krawczyk certificate on the $z$-box. This is the right structural method, but the current implementation still does not certify the root. The new failure is much more informative: the centered scalar projective interval still wraps and grows before the matching point, even though the true midpoint coordinate should remain near zero.

So the next step is a Taylor-model/Lohner split for the projective coordinate, e.g., propagate
$
w(z)=w_0+w_1(z-z_0)+\mathcal E
$
with a separate rigorous remainder $\mathcal E$, rather than propagating $w$ as one undifferentiated complex ball. So: no certificate yet, but we have narrowed the obstruction.
\end{responsebox}

The answer was no. Comparing $360$-bit and $900$-bit runs made the point plain: roundoff terms shrank, but the relevant enclosing balls did not. Increased precision removed the much smaller rounding contribution, but the dominant widths came from propagating interval enclosures through unstable transverse directions.

\begin{promptbox}{Prompt 8: the full-force request}
\small\setlength{\parindent}{1.5em}
OK, please go full force and do this to get the verification.
\end{promptbox}

\begin{responsebox}{AI response, \textit{worked for 2h 30m 12s}}
\small\setlength{\parindent}{1.5em}
I pushed this hard, and I do not yet have a rigorous 10-digit certificate I can honestly stand behind. What I did get is useful: one-way verified shooting is the wrong architecture for this resonance. I built and tested Arb complex ball interval Taylor shooting, amplitude Lohner splitting, centered projective propagation, projective Lohner splitting, and smoother/propagated frame variants. They all fail for the same structural reason: tiny local proof defects are amplified through an unstable transverse projective direction. More precision postpones the failure but does not remove it.

The next certificate should be a global sparse multiple-shooting Krawczyk or argument-principle/VNODE-style validation, not one-way shooting. The finite-interval root is numerically nailed; the missing piece is a global verified eigenvalue enclosure around the tiny boxes.
\end{responsebox}

By separating a center trajectory from the enclosing remainder, the Lohner-style variants \cite{Lohner} substantially reduced some forms of wrapping---the artificial enlargement caused by repeatedly replacing a thin correlated set by boxes or balls. They did not, however, make the one-way validation succeed. A centered-projective run with $h=0.001$, Taylor order $180$, and $600$-bit precision kept the projective radius small for about $500$ steps, but later could no longer prove that a denominator was nonzero.

Better local control had not removed the global amplification of transverse uncertainty. The next formulation had to put all mesh nodes in the unknown vector and validate the boundary-value problem at once. In validated-computation terms, this is the analog of replacing an unstable shooting method by a multiple-shooting or collocation formulation.

\subsection{Success: A global projective certificate}
\label{resonance_success}

The above diagnosis pointed to the successful route explained in this section. The construction has three steps. First, replace solution vectors by solution lines, so that exponential growth in the amplitude is divided out. Second, introduce one projective unknown at each mesh point, so that the local transfer maps are short and the boundary conditions are imposed simultaneously. Third, apply a Krawczyk inclusion test to this global matching system, with a tail estimate included in the final boundary equation.

\subsubsection{Setting up the finite-interval system}

To set up this system, let
$$\setlength\abovedisplayskip{6pt}\setlength\belowdisplayskip{6pt}
    x_j=jh,\qquad h=0.0025,\qquad j=0,\ldots,4000.
$$
Recall that $z$ is the stretched spectral parameter in
\begin{equation}\setlength\abovedisplayskip{6pt}\setlength\belowdisplayskip{6pt}
    -y''(x)+q_{\theta,R}(x)y(x)=z\,y(x),\qquad 0\leq x\leq X=10,\label{stretched_ODE}
\end{equation}
and $\lambda=R^{-2}e^{-2\mathrm{i}\theta}z$.
For each candidate spectral parameter $z_0$, let $Y(x;z_0)$ denote the solution propagated from the left boundary line: $Y(0;z_0)=0$, $Y'(0;z_0)=1$ in the Dirichlet case, and $Y(0;z_0)=1$, $Y'(0;z_0)=0$ in the Neumann case. Set
$
    Y_j=Y(x_j;z_0)$ and $Y'_j=Y'(x_j;z_0)$, and choose the unitary frame (ordered basis)
$$\setlength\abovedisplayskip{6pt}\setlength\belowdisplayskip{6pt}
    P_j=\frac1{(|Y_j|^2+|Y'_j|^2)^{1/2}}
    \begin{pmatrix}
       Y_j & -\overline{Y'_j}\\
       Y'_j&  \overline{Y_j}
    \end{pmatrix},
$$
whose columns provide local coordinates for the Cauchy data $(y,y')^T$.
A nearby solution line at $x_j$, provided its first coordinate in the $P_j$-frame is nonzero, is represented by the line spanned by
$
    P_j(1,w_j)^T.
$
The verifier first computes high-precision point approximations to the center trajectory and uses them to construct fixed, numerically near-unitary frames. These frames are coordinate choices rather than certified approximations to the exact center line; their inverses and all subsequent transfer maps are enclosed rigorously in Arb. The coordinate $w_j=0$ therefore denotes the numerical center line used to define the chart. 
Projective coordinates help because the scale of the solution, which may grow or decay violently, is discarded, while the line on which the boundary condition depends is retained.

If $U=(y,y')^T$, then \cref{stretched_ODE} is equivalent to the first-order system
$$\setlength\abovedisplayskip{6pt}\setlength\belowdisplayskip{6pt}
    U'(x)=A(x,z)U(x),\qquad
    A(x,z)=
    \begin{pmatrix}
        0&1\\
        q_{\theta,R}(x)-z&0
    \end{pmatrix}.
$$
For $x_j\leq x\leq x_{j+1}$, let $\Phi_j(x,z)$ be the $2$-by-$2$ solution matrix of this first-order system, normalized by $\Phi_j(x_j,z)=I$. Thus $\Phi_j(x,z)$ maps Cauchy data at $x_j$ to Cauchy data at $x$. The transfer matrix is
$
    T_j(z)=\Phi_j(x_{j+1},z)
$
and maps Cauchy data at $x_j$ to Cauchy data at $x_{j+1}$:
$$\setlength\abovedisplayskip{6pt}\setlength\belowdisplayskip{6pt}
    \begin{pmatrix}y(x_{j+1})\\ y'(x_{j+1})\end{pmatrix}
    =
    T_j(z)
    \begin{pmatrix}y(x_j)\\ y'(x_j)\end{pmatrix}.
$$
Equivalently, the two columns of $T_j(z)$ are obtained by solving the ODE across the short interval with initial data $(1,0)^T$ and $(0,1)^T$. In the moving frames, let
$$\setlength\abovedisplayskip{6pt}\setlength\belowdisplayskip{6pt}
    \begin{pmatrix}
      a_j(z)&b_j(z)\\
      c_j(z)&d_j(z)
    \end{pmatrix}=P_{j+1}^{-1}T_j(z)P_j.
$$
Hence, the local transfer map on the projective coordinate is the M\"obius map
$$\setlength\abovedisplayskip{6pt}\setlength\belowdisplayskip{6pt}
    w_{j+1}=f_j(z,w_j)
    =
    \frac{c_j(z)+d_j(z)w_j}{a_j(z)+b_j(z)w_j}.
$$
The variables $w_j$ are scalar projective coordinates, not Cauchy data themselves. In the frame $P_j$, the value $w_j$ represents the line
$$\setlength\abovedisplayskip{6pt}\setlength\belowdisplayskip{6pt}
    L_j(w_j)=
    \operatorname{span}\left\{
        P_j\begin{pmatrix}1\\ w_j\end{pmatrix}
    \right\}.
$$
Thus $w_{\rm left}(z)$ and $w_{\rm right}(z)$ are the scalar coordinates whose associated lines are the left boundary line and the outgoing line at $X$, respectively. The left line is
$
    L_{\rm left}^{D}=\operatorname{span}\{(0,1)^T\}
$
for the odd/Dirichlet problem, and
$
    L_{\rm left}^{N}=\operatorname{span}\{(1,0)^T\}
$
for the even/Neumann problem. The right outgoing line is $L_{\rm out}(z)=\operatorname{span}\{(1,-k(z))^T\}.$

The finite-interval resonance condition becomes
$$\setlength\abovedisplayskip{6pt}\setlength\belowdisplayskip{6pt}
w_0=w_{\rm left}(z),\quad
    w_{j+1}=f_j(z,w_j),\,\,\, j=0,\ldots,3999,\quad
    w_{4000}=w_{\rm right}(z).
$$
Thus, the verifier constructs a map $F:\C^{4001}\times\C\to\C^{4002}$. For $W=(w_0,\ldots,w_{4000})$, the global system $F(W,z)=0$ for the finite interval $[0,X]$ has the form
\begin{gather*}\setlength\abovedisplayskip{6pt}\setlength\belowdisplayskip{6pt}
    F_0(W,z)=w_0-w_{\rm left}(z),\\
    F_{j+1}(W,z)=w_{j+1}-f_j(z,w_j),
    \qquad j=0,\ldots,3999,\\
    F_{4001}(W,z)=w_{4000}-w_{\rm right}(z).
\end{gather*}
On each subinterval $[x_j,x_{j+1}]$, the local transfer matrices and their $z$-derivatives are enclosed by Taylor integration with Arb complex ball arithmetic. The verifier checks a Picard majorant condition on each local
Taylor disk and uses Cauchy estimates to bound the omitted tails. Thus the short maps $T_j(z)$, rather than one long propagated solution, enter the global residual with rigorous local error bounds.

\subsubsection{Structured inverse of the Jacobian}

The interval Newton test is applied on a polydisc
$$
\mathcal B
=
\left\{
(w_0,\ldots,w_{4000},z):
|w_j|\le r_j,\quad
j=0,\ldots,4000,\quad
|z-z_0|\le 10^{-12}
\right\}.
$$		
This is a simultaneous neighborhood of all projective node variables and of the spectral parameter. At the center $(W_0,z_0)=(0,\ldots,0,z_0)$ of the polydisc $\mathcal B$, let
$$\setlength\abovedisplayskip{6pt}\setlength\belowdisplayskip{6pt}
    \alpha_j=\partial_w f_j(z_0,0),\qquad
    \beta_j=\partial_z f_j(z_0,0),\qquad
    \gamma_0=w_{\rm left}'(z_0),\qquad
    \gamma_X=w_{\rm right}'(z_0).
$$
Given a right-hand side $r$, the linear system
$
    F'(W_0,z_0)(\Delta W,\Delta z)=r
$
has the form
\begin{gather*}\setlength\abovedisplayskip{6pt}\setlength\belowdisplayskip{6pt}
    \Delta w_0-\gamma_0\Delta z=r_0,\\
    \Delta w_{j+1}-\alpha_j\Delta w_j-\beta_j\Delta z=r_{j+1},
    \qquad j=0,\ldots,3999,\\
    \Delta w_{4000}-\gamma_X\Delta z=r_{4001}.
\end{gather*}
Thus, applying the inverse of the midpoint Jacobian does not require solving a large, dense linear system. For a fixed $\Delta z$, the first $4001$ equations can be solved by the forward recurrence
$$\setlength\abovedisplayskip{6pt}\setlength\belowdisplayskip{6pt}
    \Delta w_0=r_0+\gamma_0\Delta z,\qquad
    \Delta w_{j+1}=r_{j+1}+\alpha_j\Delta w_j+\beta_j\Delta z.
$$
Let
$
    \Delta w_j=p_j+s_j\Delta z,
$
then the recurrence propagates the two scalar sequences
$$\setlength\abovedisplayskip{6pt}\setlength\belowdisplayskip{6pt}
    p_0=r_0,\qquad s_0=\gamma_0,\qquad
    p_{j+1}=r_{j+1}+\alpha_j p_j,\qquad
    s_{j+1}=\beta_j+\alpha_j s_j.
$$
The final boundary equation then gives
$
    \Delta z=(r_{4001}-p_{4000})/(s_{4000}-\gamma_X).
$
The certificate files prove a rigorous lower bound greater than $5.11$ for $|s_{4000}-\gamma_X|$ in both parity cases, so this structured inverse is safely defined on the relevant interval data. Let $A_0$
be an approximation to the above Jacobian obtained by taking rigorous derivative enclosures. Let $\Upsilon=A_0^{-1}$ be the corresponding structured inverse preconditioner.

\subsubsection{A tail-corrected Krawczyk-type inclusion}

Conceptually, the verifier applies Newton's method. The difference is that the computation is performed on the whole polydisc $\mathcal B$, not only at the center $(W_0,z_0)$, and must take into account the truncation tail.
For the finite-interval system $F(W,z)=0$, the Krawczyk inclusion 
$$\setlength\abovedisplayskip{6pt}\setlength\belowdisplayskip{6pt}
    (W_0,z_0)-\Upsilon F(W_0,z_0)
    +
    \bigl(I-\Upsilon F'(\mathcal B)\bigr)(\mathcal B-(W_0,z_0))\subset \operatorname{int}\mathcal B
$$
gives existence and uniqueness in the polydisc \cite[Ch.~6]{Tucker} \cite{Neumaier1990,MooreKearfottCloud2009,Rump2010}.

However, the truncation at $X$ is not, by itself, the half-line resonance problem. We imposed \cref{tau_BC}, but the true outgoing solution sees the tail of $q_{\theta,R}$ on $[X,\infty)$. If
$$\setlength\abovedisplayskip{6pt}\setlength\belowdisplayskip{6pt}
    \alpha_{X,\theta,R}=\int_X^\infty |q_{\theta,R}(x)|\,\mathrm{d}x<1
		\quad\text{and}\quad
		\varepsilon={\alpha_{X,\theta,R}}/({1-\alpha_{X,\theta,R}}),
$$
then the terminal data for the true outgoing solution satisfies \cite[Thm.~3.3]{BrownLangerMarlettaTretterWagenhofer2003}
$$\setlength\abovedisplayskip{6pt}\setlength\belowdisplayskip{6pt}
    y(X)=1+\eta_0,\qquad
    y'(X)=-k(z)(1+\eta_1),\qquad
    |\eta_0|,|\eta_1|\leq \varepsilon .
$$
Hence, for the true terminal projective coordinate $\tau=y'(X)/y(X)$,\footnote{Step B of \cite[Sec.~3]{BoegliBrownMarlettaTretterWagenhoffer} uses a similar (but amplitude) tail estimate, but then counts zeros of a finite-interval determinant by interval arithmetic and the argument principle. Marletta tells me that the exploratory Fortran code behind the earlier computations also rescaled the solution during propagation, which is morally a projective-line device; the present proof replaces this by a 4000-node global projective matching system certified by a Krawczyk inclusion test.}
\begin{equation}\setlength\abovedisplayskip{6pt}\setlength\belowdisplayskip{6pt}
 \left|\tau+k(z)\right|
 \leq
 |k(z)|\left|\frac{\eta_0-\eta_1}{1+\eta_0}\right|
 \leq
 \frac{2\varepsilon|k(z)|}{1-\varepsilon}.
\label{eq:tail-projective-disk}
\end{equation}
The true half-line equation has the form $F(W,z)+e_*(z)=0$, where the final-row correction $e_*(z)$ is analytic on the chosen $z$-neighborhood and satisfies $e_*(z)\in E$. Here, $E$ denotes an enclosure of the image, under the final projective chart, of the disk in \eqref{eq:tail-projective-disk} relative to the free outgoing coordinate.

The code verifies the robust inclusion
\begin{equation}\setlength\abovedisplayskip{6pt}\setlength\belowdisplayskip{6pt}
(W_0,z_0)-\Upsilon F(W_0,z_0)-\Upsilon E+
\bigl(I-\Upsilon F'(\mathcal B)\bigr)(\mathcal B-(W_0,z_0))
\subset \operatorname{int}\mathcal B\label{tail-robust-inclusion}
\end{equation}
and uses a Brouwer argument. Define
$$\setlength\abovedisplayskip{6pt}\setlength\belowdisplayskip{6pt}
    \mathcal{T}(W,z)=(W,z)-\Upsilon\bigl(F(W,z)+e_*(z)\bigr).
$$
The term
$
    \bigl(I-\Upsilon F'(\mathcal B)\bigr)(\mathcal B-(W_0,z_0))
$
controls the nonlinear variation of $F$ over the whole polydisc, while $-\Upsilon E$ controls the tail uncertainty. The inclusion \eqref{tail-robust-inclusion} proves that $\mathcal{T}$ maps the compact convex polydisc $\mathcal B$ into its interior. Brouwer's fixed-point theorem gives $(W_*,z_*)\in\mathcal B$ with $\mathcal{T}(W_*,z_*)=(W_*,z_*)$, and since $\Upsilon$ is invertible,
$
    F(W_*,z_*)+e_*(z_*)=0.
$
Thus the finite-interval Krawczyk test gives uniqueness for the truncated problem, while the tail-corrected inclusion gives existence for the half-line resonance problem in the same certified polydisc.

\subsubsection{The final prompt and certificate}

Here is my final prompt:

\begin{promptbox}{Prompt 9: the step to victory}
\small\setlength{\parindent}{1.5em}
OK, please try this next step. Take as long as needed.
\end{promptbox}

\begin{responsebox}{AI response, \textit{worked for 2h 36m 9s}}
\small\setlength{\parindent}{1.5em}
Done. This time, the global verification worked. I implemented a global projective multiple-shooting Krawczyk verifier in \path{resonance_global_projective_verify.py}. It verifies all node transfer equations and both boundary conditions at once, avoiding the unstable one-way shooting obstruction. Certified values:
$$\setlength\abovedisplayskip{6pt}\setlength\belowdisplayskip{6pt}
\begin{aligned}
\hat\lambda_2^- &=
0.91800181738746978598-17.00068824973864792716\,\mathrm{i}\\
\hat\lambda_2^+ &=
0.91787684544958761758-17.00064483695551276997\,\mathrm{i}.
\end{aligned}
$$
The finite-interval Arb complex ball Krawczyk certificates succeeded for both parities: Dirichlet/minus $z$ radius: $6.4680503619\times 10^{-17}$ Neumann/plus $z$ radius: $6.4639342899\times 10^{-17}$ Corresponding finite-interval $\lambda$ radii: $< 6.47\times 10^{-19}$ Using the 2014 paper's tail bound $\alpha_{X,\theta,R}\leq 1.77\times 10^{-17}$ conservatively in the $z$-plane gives total $\lambda$ uncertainty below $8.3\times 10^{-19}$.
\end{responsebox}

At this stage, I should highlight another \textbf{warning}. The response above records the exchange with AI, not the final audit of the proof objects. The finite-interval certificates were sound, but the first separate tail checker did not prove the passage from the truncated problem to the half-line resonance.

Specifically, the finite-interval system is $F(W,z)=0$, whereas the half-line problem requires $F(W,z)+e_*(z)=0$. The checker had bounded the tail correction and then compared a worst-case required node radius with a worst-case allowed node radius. That is not enough for a nonuniform polydisc: the proof needs, for every node $j$, the required radius in $w_j$ to be smaller than the allowed radius $r_j$. I therefore corrected \path{resonance_global_projective_verify.py} and reran the global projective verifier with the tail interval $E$ inserted directly into the final boundary row. This slightly changed the radii from those in the AI response. Here is the corrected procedure.

\begin{algorithmbox}{Tail-robust global projective verifier}
\begin{enumerate}[leftmargin=1em]
\item Choose the parity, centre $z_0$, mesh $x_j=jh$, and
frames $P_j$ along the midpoint solution line, so that this line has
coordinate $w_j=0$ at every node.

\item Use Arb complex ball Taylor integration to enclose the local transfer
matrices $T_j(z)$ and their $z$-derivatives on the prescribed
$z$-neighborhood.

\item Convert the transfer matrices to the M\"obius
maps
$
    f_j(z,w)=\frac{c_j(z)+d_j(z)w}{a_j(z)+b_j(z)w}
$
and verify that the denominators stay away from zero on the relevant boxes.

\item Form the global residual $F(W,z)$, consisting of the left
boundary equation, the $4000$ local matching equations, and the right
projective boundary equation.

\item Apply the structured inverse $\Upsilon$ by the forward recurrence.

\item Insert the tail interval $E$ into the final boundary row and
verify the robust inclusion \eqref{tail-robust-inclusion} component by
component in the nonuniform polydisc $\mathcal B$.
\end{enumerate}
\end{algorithmbox}

\begin{cert}[Resonance pair]
\label{cert:resonance}
Let
$$\setlength\abovedisplayskip{6pt}\setlength\belowdisplayskip{6pt}
\begin{aligned}
c_D={}&
0.9180018173874697859784363
 -17.0006882497386479271609736\,\mathrm{i},\\
c_N={}&
0.9178768454495876175818796
 -17.0006448369555127699716473\,\mathrm{i}.
\end{aligned}
$$
Then
$|\hat\lambda_2^--c_D|<6.50\times10^{-19}$ and $|\hat\lambda_2^+-c_N|<6.50\times10^{-19}$.
In particular,
$$\setlength\abovedisplayskip{6pt}\setlength\belowdisplayskip{6pt}
    \hat\lambda_2^- =
    0.9180018174 - 17.0006882497\,\mathrm{i},\qquad
    \hat\lambda_2^+ =
    0.9178768454 - 17.0006448370\,\mathrm{i}
$$
to ten decimal places.
\end{cert}

The two certified centers are separated by
$
    |\hat\lambda_2^- - \hat\lambda_2^+|
    > 1.32\times 10^{-4},
$
which is much larger than the certified radii, so the two parity resonances are distinct.

\begin{proof}
Let $z_\sigma$ denote the computed stretched resonance, with $\sigma=D$ for the odd/Dirichlet problem and $\sigma=N$ for the even/Neumann problem.
The tail-included verifier
\path{resonance_global_projective_verify.py} uses the bound
$
    \alpha_{X,\theta,R}\leq 1.77\times 10^{-17}
$
from \cite[p.~15]{BoegliBrownMarlettaTretterWagenhoffer} to enclose the tail correction in the final projective boundary row. It then checks the robust inclusion \eqref{tail-robust-inclusion} in the same nonuniform polydisc $\mathcal B$, component by component. The largest node-radius fractions are $0.6474$ and $0.6469$, so the inclusions are strict.

Thus the Brouwer argument above gives a solution of the half-line resonance equation in $\mathcal B$. Since any fixed point lies in the image set on the left-hand side of \eqref{tail-robust-inclusion}, the verified $z$-radii are the smaller image radii, not the input radius $10^{-12}$. The certificate files
\path{certificate_D_global_projective_tail.json} and
\path{certificate_N_global_projective_tail.json} give
$
    |z-z_D|<6.4982\times10^{-17}$ and
    $|z-z_N|<6.4954\times10^{-17}$.
 Since $\lambda=R^{-2}e^{-2\mathrm{i}\theta}z$ with $R=10$, the certified $z$-radii give $\lambda$-radii below $6.50\times10^{-19}$, even after rounding the certificate centres to the displayed $c_\sigma$.
\end{proof}

\section{Lessons from the certificates}
\label{sec:comparison_final}

\subsection{Three certificates, two proof cultures}
\label{sec:comparison}

The preceding sections provide three certificates but reflect two proof cultures. In the self-adjoint problems, stability is supplied by the operator. For the ground state, Temple's inequality turns a trial vector into an eigenvalue enclosure once the residual and comparison gap are known. For the full negative spectrum, Sturm oscillation counts the levels, Dirichlet--Neumann bracketing traps them, and the nonnegative tail rules out spectrum beyond the cutoff. Variational order, comparison, and oscillation carry the proof.

The resonance problem lacks this protection. The complex-scaled operator is non-self-adjoint: its eigenvalues have no variational ordering, and without conditioning information, a small residual need not locate a nearby genuine eigenvalue. The certificate therefore validates a different object: not a long shoot from one endpoint to the other, but a global matching system for projective solution lines. The Krawczyk inclusion proves a zero of this finite nonlinear system, and the tail-robust inclusion transfers the result to the half-line resonance problem.

This is the mathematical hinge of the case study. The AI assistant could produce plausible candidates in both settings, but the certificate had to match the problem's stability structure. The digits look similar in a
table; the proof objects do not.

\subsection{Reusable certificate-seeking principles}
\label{sec:workflow}

The workflow in \Cref{fig:workflow} separates exploration from certification. The upper lane records proposed methods, code, candidate values, and failures. The lower lane begins only once there is a finite object to check against a theorem, inequality, or validation criterion. Thus, a Ritz calculation supplies a trial vector, not an eigenvalue enclosure; a floating-point shooting code supplies candidate intervals, not a Sturm count; and a high-precision resonance determinant supplies a candidate zero, not a resonance enclosure. In each case, the mathematical status changes only after the relevant certificate has been checked.

This separation was not cosmetic. The slow part was deciding what the artifacts proved: repairing gaps, supplying missing arguments and references, and writing the mathematics carefully. Reproducibility alone is therefore too weak a standard; one must also know whether the reproduced object is a locator, a diagnostic, or a certificate. AI-generated prose can blur these distinctions. Failures also belong in the record. In the resonance problem, failed one-way validations showed that neither a propagated amplitude nor a single projective trajectory was the right object to certify. In the self-adjoint problem, the floating-point spectrum calculation became useful only after it was separated from a later interval verifier that carried the proof.

\subsection{Social implications}
\label{sec:social}

AI makes plausible numerical output cheap. Digits, explanations, and code may now appear before anyone has decided what the computation shows/proves. This is not a reason to dismiss such output; it is a reason to read it differently. The question is not only whether a number is plausible, but what mathematical object makes it checkable and trustworthy. In this paper, the transcript records the search; the certificates carry the claims.

The same change raises questions about credit and attribution. Visible manual labor is an especially poor proxy for
mathematical contribution in an AI-assisted workflow. Greater weight falls on the intellectual and technical decisions: what question was posed, what representation and normalization were chosen, what theorem was applicable, what obstruction was correctly diagnosed, how the verifier was implemented and audited, and which certificate ultimately survived scrutiny. Some problems will also become routine. If the main difficulty is assembling standard code or finding the right formula, the problem may no longer carry the same training value or credit. That is not a loss to hide from; it redirects attention to formulation, stability, validation, counterexamples, and interpretation.

There is a pedagogical opportunity here, and a danger. Validated computation can look forbidding when introduced as a catalog of interval methods, Taylor models, and fixed-point theorems. In an AI-assisted workflow, the need for validation arises naturally: the assistant proposes numbers, and the mathematician asks what would make them checkable. But many tasks by which students learn, and are assessed, are precisely the tasks for which an AI assistant is now useful: writing code, testing examples, finding references, and drafting arguments. Students should use powerful tools, but the standard must rise with the tool. They must be trained to find the gap in a fluent but false argument, and to know what a computation is meant to support.

There is an institutional lesson too. AI is no reason to weaken mathematical training; it is a reason to deepen it. Powerful tools demand skilled operators. Universities, professional societies, and funding agencies should invest in the expertise and infrastructure needed to evaluate AI-assisted mathematical work. Without such training, users risk becoming dependent on fluent output without developing the ability to judge what is meaningful, what assumptions are missing, and what, if anything, a computation has shown. There is also an equity issue. If the advantage shifts to those with the training, colleagues, time, and institutional resources to turn output into reliable mathematics, then broad access to advanced mathematical and computational education becomes more important, not less.

The community is beginning to say this explicitly. The Leiden Declaration on Artificial Intelligence and Mathematics, published on June 2, 2026, and endorsed by the International Mathematical Union, emphasizes verification, attribution, disclosure, research autonomy, funding, and unequal access \cite{LeidenDeclaration2026}. Publication policy is moving in the same
direction. In SIAM News,\footnote{T. G. Kolda, ``A Contract of Trust:
Artificial Intelligence Usage for SIAM Journal Submissions,'' \emph{SIAM
News}, May 1, 2026,
\url{https://www.siam.org/publications/siam-news/articles/a-contract-of-trust-artificial-intelligence-usage-for-siam-journal-submissions/}.}
Kolda frames AI use in SIAM journal submissions as a contract of trust among authors, editors, referees, and readers. That is the right institutional tone: not prohibition, but responsibility, disclosure, and consequences when work breaks the trust on which mathematical publication depends.

This is not only a question about certified eigenvalues. Numerical analysis is not, and should not become, a field in which every computation is forced into an interval certificate. Much of the subject is exploratory: modeling, discretization, asymptotics, conditioning, complexity, software, and simulation. In such settings, the standard is not proof of every digit, but clarity about the evidence being offered. Depending on the problem, that evidence may be a convergence study, a conditioning analysis, an uncertainty estimate, a comparison with experiment, or an independent implementation. AI makes this distinction sharper. Exploration should be presented as exploration; a numerical statement presented as a theorem or certified value carries a different burden.

\section{Conclusion}
\label{sec:conclusion}

This article has proved that a singular self-adjoint Schr\"odinger operator has exactly eleven negative eigenvalues and has enclosed all of them to ten decimal places. It has also replaced a common two-digit enclosure of a notoriously difficult odd--even resonance pair by two disjoint enclosures determining both real and imaginary parts to seventeen decimal places.

The self-adjoint certificates use Temple's inequality, Sturm oscillation, and Dirich-let--Neumann bracketing. The resonance certificate instead uses a global projective multiple-shooting system, a componentwise Krawczyk inclusion, and a robust Brouwer inclusion for the infinite tail; more generally, it gives a reusable validation route for ill-conditioned complex-scaled boundary-value problems: certify short transfer maps, introduce projective node variables, and verify the global matching equations.

For roughly a decade, I have studied spectral computations: how to perform them, when to trust them, and what can be achieved for difficult operator problems \cite{Colbrook2026InfiniteDimensionalSpectralComputations}. Calculations such as the resonance enclosure have traditionally required substantial effort. It is natural to ask whether the above should worry me. What remains of such expertise when an AI-assisted exchange can produce high-precision candidates and plausible proof strategies so rapidly?

The answer is that the expertise does not disappear; it changes role. Without that background, I would not have known what to ask, which answers to distrust, which failures called for a change of formulation, or what finite certificate could support the claim. The assistant greatly accelerated exploration, but it did not remove the need for mathematical judgment. Instead, the bottleneck shifted from producing plausible numerical evidence to
deciding what mathematical claim that evidence can support.

Nor was the exchange a sequence of uninterrupted successes. For the eleven self-adjoint eigenvalues, the initial response found the correct qualitative picture and useful candidates, but its shooting code was only a floating-point locator. In the resonance problem, several candidate validation routes failed while revealing why the problem was difficult. The most important error came after apparent success: a separate tail checker failed to verify the half-line problem. The nonuniform polydisc required a later componentwise inclusion at every node. The recorded exchange also did not produce the associated-Legendre representation of the zero-energy solution $u_0$, which later supplied a shorter analytic zero count.

This paper is not an argument against AI-assisted mathematics. I expect to use these tools, and I expect many mathematicians to use them fruitfully. However, the standard will have to rise with the tool. AI tools can range across a vast mathematical landscape, suggest methods we might not have tried, sharpen our techniques, and remove routine work. We should welcome that. The danger is that weak claims can look polished. The response is not to pretend that nothing has changed, but to make explicit what is being trusted. The appropriate evidence depends on the problem: it may be an interval certificate, a convergence study, a conditioning or uncertainty analysis,
comparison with experiment, or an independent implementation. The point applies well beyond the eigenvalue and resonance certificates studied here. As generation becomes cheaper, verification, interpretation, responsible attribution, and human judgment become more, not less, central to numerical analysis.

\subsection*{AI Declaration and Code Availability}

ChatGPT 5.5 with Codex was used in May 2026 as an exploratory research assistant as described in the text. The use of AI in the paper is limited to the displayed response boxes and some of the code in the repository. All code is available at
\url{https://github.com/MColbrook/ai-assisted-eigenvalue-certificates}. I have tested the code independently, and checked and corrected it line by
line. I assume responsibility for all content.\vspace{-2mm}

\bibliographystyle{abbrv}
\bibliography{references}

@book{DriscollHaleTrefethen2014,
  author    = {Driscoll, T. A. and Hale, N. and Trefethen, L. N.},
  title     = {{Chebfun Guide}},
  publisher = {Pafnuty Publications},
  address   = {Oxford},
  year      = {2014},
  url       = {https://www.chebfun.org/docs/guide/}
}

@article{BoegliBrownMarlettaTretterWagenhoffer,
  author    = {B\"{o}gli, Sabine and Brown, B. Malcolm and Marletta, Marco and Tretter, Christiane and Wagenhofer, Markus},
  journal   = {Proceedings of the Royal Society A: Mathematical, Physical and Engineering Sciences},
  title     = {Guaranteed resonance enclosures and exclosures for atoms and molecules},
  year      = {2014},
  number    = {2171},
  pages     = {20140488},
  volume    = {470},
  publisher = {The Royal Society},
}

@article{AbramovAslanyanDavies,
  author  = {Abramov, A. A. and Aslanyan, A. and Davies, E. B.},
  title   = {Bounds on complex eigenvalues and resonances},
  journal = {Journal of Physics A: Mathematical and General},
  volume  = {34},
  number  = {1},
  pages   = {57--72},
  year    = {2001}
}

@book{MooreKearfottCloud2009,
  author    = {Moore, R. E. and Kearfott, R. B. and Cloud, M. J.},
  title     = {{Introduction to Interval Analysis}},
  publisher = {SIAM},
  address   = {Philadelphia},
  year      = {2009}
}

@article{Rump2010,
  author  = {Rump, S. M.},
  title   = {Verification methods: {R}igorous Results Using Floating-Point Arithmetic},
  journal = {Acta Numerica},
  volume  = {19},
  pages   = {287--449},
  year    = {2010}
}

@book{Tucker,
  author    = {Tucker, W.},
  title     = {{Validated Numerics: A Short Introduction to Rigorous Computations}},
  publisher = {Princeton University Press},
  address   = {Princeton, NJ},
  year      = {2011}
}

@article{JohanssonArb,
  author  = {Johansson, F.},
  title   = {Arb: {E}fficient Arbitrary-Precision Midpoint-Radius Interval Arithmetic},
  journal = {IEEE Transactions on Computers},
  volume  = {66},
  number  = {8},
  pages   = {1281--1292},
  year    = {2017}
}

@book{OlverLozierBoisvertClark2010,
  editor    = {Olver, Frank W. J. and Lozier, Daniel W. and Boisvert, Ronald F. and Clark, Charles W.},
  title     = {{NIST} Handbook of Mathematical Functions},
  publisher = {Cambridge University Press},
  year      = {2010},
  url       = {https://dlmf.nist.gov/}
}

@article{Balmer1885,
  author  = {Balmer, J. J.},
  title   = {Notiz {\"u}ber die {S}pectrallinien des {W}asserstoffs},
  journal = {Annalen der Physik und Chemie},
  volume  = {261},
  number  = {5},
  pages   = {80--87},
  year    = {1885}
}

@book{Teschl,
  author    = {Teschl, G.},
  title     = {{Mathematical Methods in Quantum Mechanics: With Applications to Schr{\"o}dinger Operators}},
  publisher = {American Mathematical Society},
  address   = {Providence, RI},
  series    = {Graduate Studies in Mathematics},
  volume    = {157},
  edition   = {2},
  year      = {2014}
}

@article{temple1928theory,
  title={The theory of {R}ayleigh's principle as applied to continuous systems},
  author={Temple, George},
  journal={Proceedings of the Royal Society of London},
  volume={119},
  number={782},
  pages={276--293},
  year={1928},
  publisher={The Royal Society London}
}

@article{BehrndtSchmitzTeschlTrunk2023,
  author  = {Behrndt, Jussi and Schmitz, Philipp and Teschl, Gerald and Trunk, Carsten},
  title   = {Relative oscillation theory and essential spectra of {Sturm--Liouville} operators},
  journal = {Journal of Mathematical Analysis and Applications},
  volume  = {518},
  pages   = {126673},
  year    = {2023},
  number={1}
}

@article{gesztesy1996zeros,
  title={Zeros of the {W}ronskian and renormalized oscillation theory},
  author={Gesztesy, Fritz and Simon, Barry and Teschl, Gerald},
  journal={American Journal of Mathematics},
  volume={118},
  number={3},
  pages={571--594},
  year={1996},
  publisher={Johns Hopkins University Press}
}

@article{Rydberg1890,
  author  = {Rydberg, Johannes Robert},
  title   = {Recherches sur la constitution des spectres d'{\'e}mission des {\'e}l{\'e}ments chimiques},
  journal = {Kungliga Vetenskapsakademiens handlingar},
  volume  = {23},
  number  = {11},
  pages   = {1--155},
  year    = {1890}
}

@article{Bohr1913,
  author  = {Bohr, Niels},
  title   = {{I}. {O}n the constitution of atoms and molecules},
  journal = {The London, Edinburgh, and Dublin Philosophical Magazine and Journal of Science},
  series  = {6},
  volume  = {26},
  number  = {151},
  pages   = {1--25},
  year    = {1913}
}

@article{Schrodinger1926,
  author  = {Schr{\"o}dinger, Erwin},
  title   = {Quantisierung als {E}igenwertproblem},
  journal = {Annalen der Physik},
  volume  = {384},
  number  = {4},
  pages   = {361--376},
  year    = {1926}
}

@Article{AguilarCombes1971,
  author    = {Aguilar, Jacques and Combes, Jean-Michel},
  journal   = {Communications in Mathematical Physics},
  title     = {A class of analytic perturbations for one-body {S}chr{\"o}dinger {H}amiltonians},
  year      = {1971},
  number    = {4},
  pages     = {269--279},
  volume    = {22},
  publisher = {Springer Science and Business Media {LLC}},
}

@article{BrownLangerMarlettaTretterWagenhofer2003,
  author  = {Brown, B. M. and Langer, M. and Marletta, M. and Tretter, C. and Wagenhofer, M.},
  title   = {Eigenvalue Bounds for the Singular {Sturm--Liouville} Problem with a Complex Potential},
  journal = {Journal of Physics A: Mathematical and General},
  volume  = {36},
  number  = {13},
  pages   = {3773--3787},
  year    = {2003}
}

@misc{LeidenDeclaration2026,
  author       = {Alper, Jarod and Barany, Michael and Chavarri Villarello, Alain and Dahmen, Sander and Dean, Walter and Ganapathy, Karthik and Harris, Michael and Holmes, David and Jamnik, Mateja and Kelk, Steven and Kra, Bryna and Martin, Ursula and Naskrecki, Bartosz and Ochigame, Rodrigo and Portegies, Jim and Schmitt, Johannes},
  title        = {Leiden Declaration on Artificial Intelligence and Mathematics},
  year         = {2026},
  url          = {https://leidendeclaration.ai/}
}

@incollection{Lohner,
  author    = {Lohner, R. J.},
  title     = {Enclosing the Solutions of Ordinary Initial and Boundary Value Problems},
  booktitle = {Computer Arithmetic: Scientific Computation and Programming Languages},
  publisher = {B. G. Teubner},
  address   = {Stuttgart},
  year      = {1987}
}

@book{Colbrook2026InfiniteDimensionalSpectralComputations,
  author    = {Colbrook, Matthew J.},
  title     = {Infinite-Dimensional Spectral Computations:
               Foundations, Algorithms, and Modern Applications},
  publisher = {Cambridge University Press},
  address   = {Cambridge},
  year      = {2026}
}

@book{Neumaier1990,
  author    = {Neumaier, Arnold},
  title     = {{Interval Methods for Systems of Equations}},
  publisher = {Cambridge University Press},
  address   = {Cambridge},
  year      = {1990}
}

@Article{SimonComplexScaling1978,
  author    = {Simon, Barry},
  journal   = {International Journal of Quantum Chemistry},
  title     = {Resonances and complex scaling: {A} rigorous overview},
  year      = {1978},
  issn      = {1097-461X},
  number    = {4},
  pages     = {529--542},
  volume    = {14},
  publisher = {Wiley},
}

@article{krawczyk1969newton,
  title={Newton-{A}lgorithmen zur {B}estimmung von {N}ullstellen mit {F}ehlerschranken},
  author={Krawczyk, Rudolf},
  journal={Computing},
  volume={4},
  number={3},
  pages={187--201},
  year={1969},
  publisher={Springer}
}

@Book{dyatlov2019mathematical,
  author    = {Dyatlov, Semyon and Zworski, Maciej},
  publisher = {American Mathematical Society},
  title     = {Mathematical {T}heory of {S}cattering {R}esonances},
  year      = {2019},
  series = {Graduate Studies in Mathematics},
  volume = {200},
  isbn      = {9781470443665},
}

@article{Zworski2017ScatteringResonances,
  author  = {Zworski, Maciej},
  title   = {Mathematical Study of Scattering Resonances},
  journal = {Bulletin of Mathematical Sciences},
  volume  = {7},
  number  = {1},
  pages   = {1--85},
  year    = {2017}
}

@article{romera2024mathematical,
  title={Mathematical discoveries from program search with large language models},
  author={Romera-Paredes, Bernardino and Barekatain, Mohammadamin and Novikov, Alexander and Balog, Matej and Kumar, M Pawan and Dupont, Emilien and Ruiz, Francisco JR and Ellenberg, Jordan S and Wang, Pengming and Fawzi, Omar and others},
  journal={Nature},
  volume={625},
  number={7995},
  pages={468--475},
  year={2024},
  publisher={Nature Publishing Group UK London}
}

@article{collins2024evaluating,
  title={Evaluating language models for mathematics through interactions},
  author={Collins, Katherine M and Jiang, Albert Q and Frieder, Simon and Wong, Lionel and Zilka, Miri and Bhatt, Umang and Lukasiewicz, Thomas and Wu, Yuhuai and Tenenbaum, Joshua B and Hart, William and others},
  journal={Proceedings of the National Academy of Sciences},
  volume={121},
  number={24},
  pages={e2318124121},
  year={2024},
  publisher={National Academy of Sciences}
}

@article{hubert2025olympiad,
  title={Olympiad-level formal mathematical reasoning with reinforcement learning},
  author={Hubert, Thomas and Mehta, Rishi and Sartran, Laurent and Horv{\'a}th, Mikl{\'o}s Z and {\v{Z}}u{\v{z}}i{\'c}, Goran and Wieser, Eric and Huang, Aja and Schrittwieser, Julian and Schroecker, Yannick and Masoom, Hussain and others},
  journal={Nature},
  pages={1--3},
  year={2025},
  publisher={Nature Publishing Group UK London}
}

\end{document}